\input amstex
\documentstyle{amsppt}
\magnification=\magstephalf \hsize = 6.5 truein \vsize = 9 truein
\vskip 3.5 in

\NoBlackBoxes
\TagsAsMath

\def\label#1{\par%
        \hangafter 1%
        \hangindent .75 in%
        \noindent%
        \hbox to .75 in{#1\hfill}%
        \ignorespaces%
        }

\newskip\sectionskipamount
\sectionskipamount = 24pt plus 8pt minus 8pt
\def\sectionskip{\vskip\sectionskipamount}
\define\sectionbreak{%
        \par  \ifdim\lastskip<\sectionskipamount
        \removelastskip  \penalty-2000  \sectionskip  \fi}
\define\section#1{%
        \sectionbreak   
        \subheading{#1}%
        \bigskip
        }

\redefine\qed{{\unskip\nobreak\hfil\penalty50\hskip2em\vadjust{}\nobreak\hfil
    $\square$\parfillskip=0pt\finalhyphendemerits=0\par}}

        
        \let    \< = \langle
        \let    \> = \rangle

\define\op#1{\operatorname{\fam=0\tenrm{#1}}} 
        \define         \a              {\alpha}
        \redefine       \b              {\beta}
        \redefine       \d              {\delta}
        \redefine       \D              {\Delta}
        \define         \e              {\varepsilon}
        \define         \E              {\op {E}}
        \define         \g              {\gamma}
        \define         \G              {\Gamma}
        \redefine       \l              {\lambda}
        \redefine       \L              {\Lambda}
        \define         \n              {\nabla}
        \redefine       \var            {\varphi}
        \define         \s              {\sigma}
        \redefine       \Sig            {\Sigma}
        \redefine       \t              {\tau}
        \define         \th             {\theta}
        \redefine       \O              {\Omega}
        \redefine       \o              {\omega}
        \define         \z              {\zeta}
        \define         \k              {\kappa}
        \redefine       \i              {\infty}
        \define         \p              {\partial}
        \define         \vsfg           {\midspace{0.1 truein}}

\topmatter
\title Counting words of minimum length in an automorphic orbit
\endtitle
\author Donghi Lee
\endauthor

\address {Department of Mathematics, Pusan National University, San-30 Jangjeon-Dong, Geumjung-Gu, Pusan, 609-735, Korea}
\endaddress

\email {donghi\@pusan.ac.kr}
\endemail


\abstract {Let $u$ be a cyclic word in a free group $F_n$ of
finite rank $n$ that has the minimum length over all cyclic words
in its automorphic orbit, and let $N(u)$ be the cardinality of the
set $\{v: |v|=|u|$ and $v= \phi(u)$ for some $\phi \in \text
{Aut}F_n\}$. In this paper, we prove that $N(u)$ is bounded by a
polynomial function with respect to $|u|$ under the hypothesis
that if two letters $x, \, y$ with $x \neq y^{\pm 1}$ occur in
$u$, then the total number of occurrences of $x^{\pm 1}$ in $u$ is
not equal to the total number of occurrences of $y^{\pm 1}$ in
$u$. A complete proof without the hypothesis would yield the
polynomial time complexity of Whitehead's algorithm for $F_n$.}
\endabstract
\endtopmatter


\document
\baselineskip=24pt

\heading 1. Introduction
\endheading

Let $F_n$ be the free group of finite rank $n$ on the set $\{x_1,
x_2, \dots, x_n\}$. We denote by $\Sig$ the set of {\it letters}
of $F_n$, that is, $\Sig = \{x_1, x_2, \dots, x_n\}^{\pm 1}$. As
in [1, 5], we define a {\it cyclic word} to be a cyclically
ordered set of letters with no pair of inverses adjacent. The {\it
length} $|w|$ of a cyclic word $w$ is the number of elements in
the cyclically ordered set. For a cyclic word $w$ in $F_n$, we
denote the automorphic orbit $\{\psi(w): \psi \in \text
{Aut}F_n\}$ by $\text {Orb}_{\text {Aut}F_n}(w)$.

The purpose of this paper is to provide a partial solution of the
following problem raised by Myasnikov--Shpilrain [6]:

\proclaim {Problem} Let $u$ be a cyclic word in $F_n$ which has
the minimum length over all cyclic words in its automorphic orbit
$\text {Orb}_{\text {Aut}F_n}(u)$, and let $N(u)$ be the
cardinality of the set $\{v \in \text {Orb}_{\text {Aut}F_n}(u):
|v|=|u|\}$. Then is $N(u)$ bounded by a polynomial function with
respect to $|u|$?
\endproclaim

This problem was settled in the affirmative for $F_2$ by
Myasnikov--Shpilrain [6], and Khan [3] improved their result by
showing that $N(u)$ has the sharp bound of $8|u|-40$ for $F_2$.
The problem was motivated by the complexity of Whitehead's
algorithm which decides whether, for given two elements in $F_n$,
there is an automorphism of $F_n$ that takes one element to the
other. Indeed, a complete positive solution to the problem would
yield that Whitehead's algorithm terminates in polynomial time
with respect to the maximum length of the two words in question
(see [6, Proposition 3.1]). Recently, Kapovich--Schupp--Shpilrain
[2] proved that Whitehead's algorithm has strongly linear time
generic-case complexity. In the present paper, we prove for $F_n$
with $n \ge 2$ that $N(u)$ is bounded by a polynomial function
with respect to $|u|$ under the following

\proclaim {Hypothesis 1.1} (i) A cyclic word $u$ has the minimum
length over all cyclic words in its automorphic orbit $\text {
Orb}_{\text { Aut}F_n}(u)$.

(ii) If two letters $x_i$ (or $x_i^{-1}$) and $x_j$ (or
$x_j^{-1}$) with $i < j$ occur in $u$, then the total number of
$x_i^{\pm 1}$ occurring in $u$ is strictly less than the total
number of $x_j^{\pm 1}$ occurring in $u$.
\endproclaim

Before we state our theorems, we would like to establish several
notation and definitions. As in [1, 5], for $A, \, B \subseteq
\Sig$, we write $A + B$ for $A \cup  B$ if $ A \cap B =
\emptyset$, and $A - B$ for $A \cap B^c$ if $ B \subseteq A$,
where $ B^c$ is the complement of $ B$ in $\Sig$. We define a {\it
Whitehead automorphism} $\s$ of $F_n$ as an automorphism of one of
the following two types (cf. [4, 7]):

\roster \item"(W1)" $\s$ permutes elements in $\Sig$. \item"(W2)"
$\s$ is defined by a set $A \subset \Sig$ and a multiplier $a \in
\Sig$ with both $a, a^{-1} \notin A$ in such a way that if $x \in
\Sig$ then (a) $\s(x)=xa$ provided $x \in A$ and $x^{-1} \notin
A$; (b) $\s(x)=a^{-1}xa$ provided both $x,\, x^{-1} \in A$; (c)
$\s(x)=x$ provided both $x,\, x^{-1} \notin A$.
\endroster
If $\s$ is of the second type, then we write $\s=(A, a)$. By
$(\bar {A}, a^{-1})$, we mean the Whitehead automorphism $(\Sig -
A - a^{\pm 1}, a^{-1})$. It is then easy to see that $(A,
a)(w)=(\bar {A}, a^{-1})(w)$ for any cyclic word $w$ in $F_n$.

For a Whitehead automorphism $\s$ of the second type, we define the degree of $\s$ as follows:

\proclaim {Definition 1.2} Let $\s=(A, a)$ be a Whitehead
automorphism of $F_n$ of the second type. Put $A'=\{i: \text
{either} \ x_i \in  A \ \text {or} \ x_i^{-1} \in  A, \ \text {but
not both} \}$. Then the {\it degree of $\s$} is defined to be
$\max {A}'$. If $A'=\emptyset$, then the {\it degree of $\s$} is
defined to be zero.
\endproclaim

Let $w$ be a fixed cyclic word in $F_n$ that satisfies Hypothesis
1.1 (i). For two letters $x, \, y \in \Sig$, we say that {\it $x$
depends on $y$ with respect to $w$} if, for every Whitehead
automorphism $(A, a)$ of $F_n$ such that
$$a \notin \{x^{\pm 1}, y^{\pm 1}\}, \ \{y^{\pm 1}\} \cap  A \neq \emptyset,
\ \text{and} \ \exists v \in \text {Orb}_{\text {Aut}F_n}(w): |(A,
a)(v)|=|v|=|w|,$$ we have $\{x^{\pm 1}\} \subseteq A$. Then we
have the following

\proclaim{Claim} If $x$ depends on $y$ with respect to $w$, then
$y$ depends on $x$ with respect to $w$.
\endproclaim

\demo{Proof} Suppose on the contrary that $y$ does not depend on
$x$. Then there exists a Whitehead automorphism $(A, a)$ of $F_n$
such that $a \notin \{x^{\pm 1}, y^{\pm 1}\}$, $x^{\pm 1} \cap A
\neq \emptyset$, $|(A, a)(v)|=|v|=|w|$ for some $v \in \text
{Orb}_{\text {Aut}F_n}(w)$, but such that $y^{\pm 1} \nsubseteq
A$. Then $|(\bar{A}, a^{-1})(v)|=|v|=|w|$ and $y^{\pm 1} \cap
\bar{A} \neq \emptyset$. Since $x$ depends on $y$, $x^{\pm 1}
\subseteq \bar{A}$. This gives $x^{\pm 1} \cap A = \emptyset$,
which is a contradiction. \qed
\enddemo

We then construct the {\it dependence graph $\G_w$ of $w$} as
follows: Take the vertex set as $\Sig$, and connect two distinct
vertices $x,\, y \in \Sig$ by a non-oriented edge if either
$y=x^{-1}$ or $y$ depends on $x$ with respect to ~$w$. Let $C_i$
be the connected component of $\G_w$ containing $x_i$. Here, we
make the following

\proclaim{Remark} (i) $\G_w=\G_v$ for any $v \in \text
{Orb}_{\text {Aut}F_n}(w)$ with $|v|=|w|$.

(ii) If $x_i$ depends on $x_j$, then $C_i=C_j$.

(iii) If $x_j^{\pm 1} \in C_i$ with $i \neq j$, then every
Whitehead automorphism $(A, a)$ such that either $x_i \in A$ or
$x_i^{-1} \in A$ but not both and such that $|(A, a)(v)|=|v|=|w|$
for some $v \in \text {Orb}_{\text {Aut}F_n}(w)$ must have the
multiplier $a$ only in $C_i$, for otherwise $x_j^{\pm 1} \subseteq
A$ but then $x_j^{\pm 1} \nsubseteq \bar{A}$, which is a
contradiction because $x_i^{\pm 1} \cap \bar{A} \neq \emptyset$.
\endproclaim

Clearly there exists a unique factorization
$$w=v_1v_2\cdots v_k \ \text {(without cancellation)},$$
where each $v_i$ is a non-empty (non-cyclic) word consisting of
letters in $C_{j_i}$ with $C_{j_i} \neq C_{j_{i+1}}$ ($i$ mod
$k$). The subword $v_i$ is called a {\it $C_{j_i}$-syllable} of
$w$. By the {\it $C_i$-syllable length of $w$} denoted by
$|w|_{C_i}$, we mean the total number of $C_i$-syllables of $w$.

For Theorem 1.4, we suppose further that a cyclic word $u$
satisfies the following

\proclaim{Hypothesis 1.3} (i) The $C_n$-syllable length
$|u|_{C_n}$ of $u$ is minimum over all cyclic words in the set
$\{v \in \text {\rm Orb}_{\text {\rm Aut}F_n}(u): |v|=|u|\}$.

(ii) If the index $j$ ($1 \le j \le n-1$) is such that $C_j \neq
C_k$ for all $k>j$, then the $C_j$-syllable length $|u|_{C_j}$ of
$u$ is minimum over all cyclic words in the set $\{v \in \text
{Orb}_{\text {Aut}F_n}(u): |v|=|u| \ \text{and} \
|v|_{C_k}=|u|_{C_k} \ \text{for all} \ k>j \}$.
\endproclaim

For an easy example, consider the cyclic words $u=x_1^2x_2^3 x_3^4
x_4^5$ and $v=x_1x_2^3x_1x_3^4 x_4^5$ in $F_4$. Clearly $v$ is an
automorphic image of $u$ with $|v|=|u|$, so $\G_u=\G_v$. The
dependence graph $\G_u=\G_v$ has four distinct connected
components, each $C_i$ of which contains only $x_i^{\pm 1}$. Then
$u$ satisfies Hypotheses 1.1 and 1.3, whereas $v$ satisfies
Hypotheses ~1.1 and 1.3 ~(i) but not Hypothesis ~1.3 ~(ii),
because the $C_1$-syllable length of $v$ can be decreased without
changing $|v|$ and $|v|_{C_i}$ for all $i>1$.

For another example, let $u=x_1^2x_2^3x_3^2x_4x_3^{-1}x_4x_3x_4^3$
and $v=x_1^2x_3^2x_2^3x_4x_3^{-1}x_4x_3x_4^3$ be cyclic words in
$F_4$. Then $v$ is an automorphic image of $u$ with $|v|=|u|$, so
$\G_u=\G_v$. In the dependence graph $\G_u=\G_v$, there are three
distinct connected components $C_1$, $C_2$, $C_3=C_4$. While $u$
satisfies Hypotheses 1.1 and 1.3, $v$ does not satisfy Hypothesis
1.3 (i), because the $C_4$-syllable length of $v$ can be decreased
without changing $|v|$.

Now we are ready to state our theorems, whose proofs will appear
in Sections 3--4.

\proclaim {Theorem 1.4} Let $u$ be a cyclic word in $F_n$ that
satisfies Hypotheses 1.1 and 1.3. Let $\s_i$, $i=1, \dots, \ell$,
be Whitehead automorphisms of the second type such that $|\s_i
\cdots \s_1(u)|=|u|$ for all $i$. Then there exist Whitehead
automorphisms $\t_1, \t_2, \dots, \t_s$ of the second type such
that
$$\s_\ell \cdots \s_2 \s_1(u) = \t_s \cdots \t_2 \t_1(u),
$$ where $\max \limits_{1 \le i \le \ell} \deg \s_i \ge
\deg \t_s \ge \deg \t_{s-1} \ge \cdots \ge \deg \t_1$, and $|\t_j
\cdots \t_1(u)|=|u|$ for all $j=1, \dots, s$.
\endproclaim

\proclaim {Theorem 1.5} Let $u$ be a cyclic word in $F_n$ that
satisfies Hypothesis 1.1, and let $N(u)$ be the cardinality of the
set $\{v \in \text {Orb}_{\text {Aut}F_n}(u): |v|=|u|\}$. Then
$N(u)$ is bounded by a polynomial function of degree $n(5n-7)/2$
with respect to $|u|$.
\endproclaim

The main idea of the present paper is to prove that the action of
an automorphism on an element which satisfies Hypotheses ~1.1 and
1.3 can be factored into a composition of automorphisms of
ascending degrees, which will be achieved through Lemmas ~3.1, 3.2
and Theorem 1.4. The proof of Theorem ~1.4 will proceed by double
induction on $\ell$ and $r$, where $\ell$ is the length of the
chain $\s_\ell \cdots \s_2 \s_1$ and $r=\max\limits_{1 \le i \le
\ell} \deg \s_i$, with Lemma ~3.1 (the case for $\ell=2$ and any
$r$) and Lemma ~3.2 (the case for $r=1$ and any $\ell$) as the
base steps of the induction.

Let $N_k(u)$ be the cardinality of the set $\{\phi (u): \phi$ can
be represented as a composition $\t_s \cdots \t_1$ ($s \in \Bbb
N$) of Whitehead automorphisms $\t_i$ of $F_n$ of degree $k$ such
that $|\t_i \cdots \t_1(u)|=|u|$ for all $i=1, \dots, s \}$. Then
bounding $N(u)$ reduces to bounding each $N_k(u)$, as will be
shown in the proof of Theorem ~1.5 using the result of Theorem
~1.4. Lemma ~4.1 will be devoted to bounding $N_0(u)$, and Lemma
~4.2 will show that $N_k(u)$ for $k \ge 1$ is at most $N_0(V_u)$,
where $V_u$ is a certain sequence of cyclic words constructed from
$u$, thus bounding $N_k(u)$ for $k \ge 1$. Furthermore in Theorem
~1.5 we will specifically give a bound for the degree of a
polynomial bounding $N(u)$.

\heading 2. Preliminaries
\endheading

We begin this section by setting some notation. Let $w$ be a fixed
cyclic word in $F_n$. As in [1], for $x, y \in \Sig$, $x. \, y$
denotes the total number of occurrences of the subwords $xy^{-1}$
and $yx^{-1}$ in $w$. For $ A , \,  B \subseteq \Sig$, $ A. \,  B$
means the sum of $a.\, b$ for all $a \in  A$, $b \in  B$. Then
obviously $a.\Sig$ is equal to the total number of $a^{\pm 1}$
occurring in $w$. For two automorphisms $\phi$ and $\psi$ of
$F_n$, by writing $\phi \equiv \psi$ we mean the equality of
$\phi$ and $\psi$ over all cyclic words in $F_n$, that is,
$\phi(v)=\psi(v)$ for every cyclic word $v$ in $F_n$.

We now establish two technical lemmas which will play a
fundamental role in the proofs in Sections 3 and 4.

\proclaim {Lemma 2.1} Let $u$ be a cyclic word in $F_n$ that
satisfies Hypothesis 1.1 (i), and let $\s=(A, a^{-1})$ and $\t=(
B, b)$ be Whitehead automorphisms of $F_n$ such that
$|\s(u)|=|\t(u)|=|u|$. Put $A=C + E$ and $B =D + E$, where $ E = A
\cap B$. Then

\noindent (i) if $a^{-1} = b$, then $|( E, a^{-1})(u)|=|u|$;

\noindent (ii) if $a^{-1} \neq b$, $a^{\pm 1} \notin  B$ and $b
\notin  A$, then $|(C, a^{-1})(u)|= |(D, b)(u)|=|u|$.
\endproclaim

\demo {Proof} It follows from [1, p.255] that
$$\cases |\s(u)|-|u|=(A+a^{-1}).(A+a^{-1})'-a.\Sig;\\
|\t(u)|-|u|=(B+b).(B+b)'-b.\Sig,
\endcases$$
where $(A +a^{-1})'=\Sig - (A+a^{-1})$ and $(B+b)'=\Sig - (B+b)$.
Since $|\s(u)|=|\t(u)|=|u|$, we have $(A+a^{-1}).(
A+a^{-1})'-a.\Sig=(B+b).(B+b)'-b.\Sig=0$, so that
$$(A+a^{-1}).(A+a^{-1})'+(B+b).(B+b)'-a.\Sig-b.\Sig=0.$$

Following the notation in [1, p.257], we write $A_1= A+a^{-1}$,
$A_2=(A+a^{-1})'$, $B_1= B+b$, $ B_2=(B+b)'$ and $P_{ij}= A_i \cap
B_j$. Then as in [1, ~p.257], we have
$$\cases P_{11} . P_{11}'+ P_{22} . P_{22}'-a.\Sig -b.\Sig=0;\\
 P_{12} . P_{12}'+ P_{21} . P_{21}'-a.\Sig -b.\Sig=0,
\endcases
\tag 2.1$$ where $P_{ij}'=\Sig - P_{ij}$.

For (i), assume that $a^{-1} = b$. Then we have $a^{-1} \in  P_{11}$ and $a \in  P_{22}$. It follows from the first equality of (2.1) that
$$\split
 P_{11} .  P_{11}'+  P_{22} .  P_{22}'-a.\Sig -a.\Sig
&=( P_{11} .  P_{11}'-a.\Sig)+( P_{22} .  P_{22}'-a.\Sig)\\
&=|( P_{11}-a^{-1}, a^{-1})(u)|-|u|+|( P_{22}-a, a)(u)|-|u|=0.
\endsplit$$
Since both $|( P_{11}-a^{-1}, a^{-1})(u)|-|u| \ge 0$ and $|(
P_{22}-a, a)(u)|-|u| \ge 0$ by Hypothesis 1.1 (i), we must have
$|( P_{11}-a^{-1}, a^{-1})(u)|=|u|$, that is, $|( E,
a^{-1})(u)|=|u|$, as required.

For (ii), assume that $a^{-1} \neq b$, $a^{\pm 1} \notin  B$ and
$b \notin  A$. Then we have $a^{-1} \in  P_{12}$, $a \notin
P_{12}$, $b \in  P_{21}$ and $b^{-1} \notin  P_{21}$. Hence the
second equality of (2.1) gives us that
$$\split
 P_{12} .  P_{12}'+  P_{21} .  P_{21}'-a.\Sig -b.\Sig
&=( P_{12} .  P_{12}'-a.\Sig)+( P_{21} .  P_{21}'-b.\Sig)\\
&=|( P_{12}-a^{-1}, a^{-1})(u)|-|u|+|( P_{21}-b, b)(u)|-|u|=0.
\endsplit$$
As above, it follows from Hypothesis 1.1 (i) that $|(
P_{12}-a^{-1}, a^{-1})(u)|=|u|$ and $|( P_{21}-b, b)(u)|=|u|$.
Since $ P_{12}-a^{-1}= C$ and $ P_{21}-b= D$, we have $|( C,
a^{-1})(u)|= |( D, b)(u)|=|u|$, as desired. \qed
\enddemo

\proclaim {Lemma 2.2} Let $u$ be a cyclic word in $F_n$ that
satisfies Hypothesis 1.1, and let $\s=( A, a)$ be a Whitehead
automorphism of $F_n$ such that $|\s (u)|=|u|$. Then $a.\Sig >
b.\Sig$ for every $b \in  A$ with $b^{-1} \notin  A$.
\endproclaim

\demo {Proof} In view of the assumption $|\s (u)|=|u|$ and [1,
p.255], we have $0=|\s(u)|-|u|=( A +a).( A+a)'-a.\Sig$, where $(
A+a)'=\Sig - (A+a)$, so that $(A +a).(A+a)'=a.\Sig$. Now let $b
\in  A$ with $b^{-1} \notin  A$. Then for the Whitehead
automorphism $\t=( A+a-b, b)$, we have $0 \le |\t(u)|-|u|=( A
+a).( A+a)'-b.\Sig$. Hence $(A +a).(A+a)'\ge b.\Sig$; thus $a.\Sig
\ge b.\Sig$. Here, the equality $a.\Sig = b.\Sig$ cannot occur by
Hypothesis ~1.1 ~(ii); therefore $a.\Sig > b.\Sig$. \qed
\enddemo

\proclaim {Remark} By Lemma 2.2, if $u$ is a cyclic word in $F_n$
that satisfies Hypothesis 1.1 and $\s=( A, a)$ is a Whitehead
automorphism of $F_n$ such that $|\s (u)|=|u|$, then $\deg \s$ is
at most $n-1$.
\endproclaim

\heading 3. Proof of Theorem 1.4
\endheading

The aim of this section is to prove Theorem ~1.4. The proof of
Theorem ~1.4 will proceed by double induction on $\ell$ and $r$,
where $\ell$ is the length of the chain $\s_\ell \cdots \s_2 \s_1$
and $r=\max\limits_{1 \le i \le \ell} \deg \s_i$. Lemma ~3.1 deals
with the case for $\ell=2$ and any $r$ as one of the base steps of
the induction. As the other base step, Lemma ~3.2 deals with the
case for $r=1$ and any $\ell$.

\proclaim {Lemma 3.1} Let $u$ be a cyclic word in $F_n$ that
satisfies Hypothesis 1.1, and let $\s_1=( A, a)$ and $\s_2=( B,
b)$ be Whitehead automorphisms of $F_n$ such that $|\s_2
\s_1(u)|=|\s_1(u)|=|u|$. Suppose that $\deg \s_1 > \deg \s_2$.
Then there exist Whitehead automorphisms $\t_1, \dots, \t_s$ of
$F_n$ of the second type such that
$$\s_2 \s_1 \equiv \t_s \cdots \t_2 \t_1,$$
where $\deg \s_1 = \deg \t_s \ge \cdots \ge \deg \t_1$ and $|\t_i \cdots \t_1(u)|=|u|$ for all $i=1, \dots, s$.
\endproclaim

\demo {Proof} It suffices to prove that there exist Whitehead
automorphisms $\g_1, \dots, \g_t$ of $F_n$ such that
$$\s_2 \s_1 \equiv \g_t \cdots \g_2 \g_1,$$
where the index $t$ is at most 3, $|\g_i \cdots \g_1(u)|=|u|$ for
all $i=1, \dots, t$, and either $\deg \s_1 = \deg \g_t > \deg
\g_j$ for all $j=1, \dots, t-1$ or otherwise $\deg \s_1 = \deg
\g_i$ for all $i=1, \dots, t$. Put $u'= \s_1 (u)$; then
$|\s_1^{-1}(u')|=|\s_2(u')|=|u|$, that is, $$|( A, a^{-1})(u')|=|(
B, b)(u')|=|u|. \tag 3.1
$$
Also put $c=x_{\deg \s_1}$. Upon replacing $( A, a)$, $( B, b)$ by
$(\bar { A}, a^{-1})$, $(\bar { B}, b^{-1})$, respectively, if
necessary, where $\bar { A}= \Sig - { A} - a^{\pm 1}$ and $\bar {
B}= \Sig - { B} - b^{\pm 1}$, we may assume that $c \in  A$ and
$c^{\pm 1} \notin  B$ (clearly $c^{-1} \notin  A$). By Lemma 2.2,
we have $a.\Sig > c.\Sig$; hence either $a^{\pm 1} \notin  B$ or
$a^{\pm 1} \in  B$, for otherwise $\deg \s_2 > \deg \s_1$,
contrary to the hypothesis $\deg \s_1 > \deg \s_2$.

We first treat four cases for $a^{\pm 1} \notin  B$ and then four
cases for $a^{\pm 1} \in  B$ according to whether $b$ or $b^{-1}$
belongs to $ A$. For convenience, we write $ A= C +  E$ and $ B =
D +  E$, where $ E =  A \cap  B$.

\proclaim {Case 1} $a^{\pm 1} \notin  B$ and $b^{\pm 1} \notin  A$.
\endproclaim

We consider two cases corresponding to whether or not $ E$ is the empty set.

\proclaim {Case 1.1} $ E = \emptyset$.
\endproclaim

\proclaim {Case 1.1.1} $a=b$.
\endproclaim

It follows from [5, relation R2] that $\s_2 \s_1 \equiv ( A +  B, a)$.

\proclaim {Case 1.1.2} $a \neq b$.
\endproclaim

By [5, relation R3], we have $\s_2 \s_1 \equiv ( A, a) ( B, b)$.

\proclaim {Case 1.2} $ E \neq \emptyset$.
\endproclaim

\proclaim {Case 1.2.1} $a=b$.
\endproclaim

In view of (3.1) and Lemma 2.1 (ii), we have $|( C,
a^{-1})(u')|=|u|$. Since $( C, a^{-1})(u')= ( E, a)(u)$, we have
$|( E, a)(u)|=|u|$; hence
$$\split
\s_2 \s_1
& \equiv ( B, a) [( C, a) ( E, a)] \equiv [( B, a) ( C, a)] ( E, a)\\
& \equiv ( C +  B, a) ( E, a) \quad \text {by Case 1.1.1},
\endsplit
$$
where $\deg \s_1=\deg ( C +  B, a)> \deg ( E, a)$.

\proclaim {Case 1.2.2} $a^{-1}=b$.
\endproclaim

Lemma 2.1 (i) together with (3.1) gives us that $|( E,
a^{-1})(u')|=|u|$, so that $|( C, a)(u)|=|u|$; thus
$$\split
\s_2 \s_1
& \equiv ( B, a^{-1}) [( E, a) ( C, a)] \equiv [( B, a^{-1}) ( E, a)] ( C, a) \equiv ( D, a^{-1}) ( C, a) \\
& \equiv ( C, a) ( D, a^{-1}) \quad \text {by Case 1.1.2},
\endsplit
$$
where $\deg \s_1=\deg ( C, a)> \deg ( D, a^{-1})$.

\proclaim {Case 1.2.3} $a^{\pm 1} \neq b$.
\endproclaim

As in Case 1.2.1, we have $|( E, a)(u)|=|u|$; hence
$$\split
\s_2 \s_1
& \equiv ( B, b) [( C, a) ( E, a)] \equiv [( B, b) ( C, a)] ( E, a) \\
& \equiv [( C, a) ( B, b)] ( E, a) \quad \text {by Case 1.1.2},
\endsplit
$$
where $\deg \s_1=\deg ( C, a)> \deg ( B, b)$, $\deg ( E, a)$.

\proclaim {Case 2} $a^{\pm 1} \notin  B$, $b \notin  A$ and $b^{-1} \in  A$.
\endproclaim

We consider this case dividing into two cases according to whether or not $ E$ is the empty set.

\proclaim {Case 2.1} $ E = \emptyset$.
\endproclaim

It follows from [5, relation R4] that $\s_2 \s_1 \equiv ( A +  B, a) ( B, b)$, where $\deg \s_1=\deg ( A +  B, a)> \deg ( B, b)$.

\proclaim {Case 2.2} $ E \neq \emptyset$.
\endproclaim

As in Case 1.2.1, we have $|( E, a)(u)|=|u|$; then
$$\split
\s_2 \s_1
& \equiv ( B, b) [( C, a) ( E, a)] \equiv [( B, b) ( C, a)] ( E, a) \\
& \equiv [( C +  B, a) ( B, b)] ( E, a) \quad \text {by Case 2.1},
\endsplit
$$
where $\deg \s_1=\deg ( C +  B, a)> \deg ( B, b)$, $\deg ( E, a)$.

\proclaim {Case 3} $a^{\pm 1} \notin  B$, $b \in  A$ and $b^{-1} \notin  A$.
\endproclaim

Since $\s_2 \s_1 \equiv ( B, b) (\bar { A}, a^{-1})$, we can apply Case 2.2 to get
$$\s_2 \s_1 \equiv ( B, b) (\bar { A}, a^{-1}) \equiv ((\bar { A} \setminus  B) +  B, a^{-1}) ( B, b) (\bar { A} \cap  B, a^{-1}).$$
Here, since $(\bar { A} \setminus  B) +  B=\Sig- C-a^{\pm 1}$ and $\bar { A} \cap  B = D$, we have
$$\s_2 \s_1 \equiv (\Sig- C-a^{\pm 1}, a^{-1}) ( B, b) ( D, a^{-1}) \equiv ( C, a) ( B, b) ( D, a^{-1}), $$
where $\deg \s_1=\deg ( C, a)> \deg ( B, b)$, $\deg ( D, a^{-1})$.

\proclaim {Case 4} $a^{\pm 1} \notin  B$ and $b^{\pm 1} \in  A$.
\endproclaim

By Case 1.2.3 applied to  $\s_2 \s_1 \equiv ( B, b) (\bar { A}, a^{-1})$, we have
$$\s_2 \s_1 \equiv ( B, b) (\bar { A}, a^{-1}) \equiv (\bar { A} \setminus  B, a^{-1}) ( B, b) (\bar { A} \cap  B, a^{-1}).$$
From the observation that $\bar { A} \setminus  B=\Sig-( C+ B)-a^{\pm 1}$ and $\bar { A} \cap  B= D$, it follows that
$$\s_2 \s_1 \equiv (\Sig-( C+ B)-a^{\pm 1}, a^{-1}) ( B, b) ( D, a^{-1}) \equiv ( C +  B, a) ( B, b) ( D, a^{-1}),$$
where $\deg \s_1=\deg ( C +  B, a)> \deg ( B, b)$, $\deg ( D, a^{-1})$.

\proclaim {Case 5} $a^{\pm 1} \in  B$ and $b^{\pm 1} \notin  A$.
\endproclaim

Since $\s_2 \s_1 \equiv (\bar { B}, b^{-1}) ( A, a)$, we have $|(
A, a^{-1})(u')|=|(\bar { B}, b^{-1})(u')|=|u|$. This implies by
Lemma ~2.1 ~(ii) that $|(\bar { B} \setminus  A,
b^{-1})(u')|=|u|$, so that
$$\s_2 \s_1 \equiv (\bar { B}, b^{-1}) ( A, a) \equiv [( A \cap \bar { B}, b^{-1}) (\bar { B} \setminus  A, b^{-1})] ( A, a).$$
Here, by Case 1.1.2, we have $(\bar { B} \setminus  A, b^{-1}) ( A, a) \equiv ( A, a)(\bar { B} \setminus  A, b^{-1})$; thus
$$\s_2 \s_1 \equiv ( A \cap \bar { B}, b^{-1}) ( A, a) (\bar { B} \setminus  A, b^{-1}).$$
Since $ A \cap \bar {B}= C$ and $\bar {B} \setminus  A=\Sig-( C+
B)-b^{\pm 1}$, we finally have
$$\s_2 \s_1 \equiv (C, b^{-1})(A, a)(C +  B, b),$$
where $\deg \s_1=\deg (C, b^{-1})=\deg (A, a)= \deg (C +  B, b)$.

\proclaim {Case 6} $a^{\pm 1} \in  B$, $b \notin  A$ and $b^{-1} \in  A$.
\endproclaim

\proclaim {Case 6.1} $c = b^{-1}$.
\endproclaim

By Case 3 applied to $\s_2 \s_1 \equiv (\bar{B}, b^{-1}) (A, a)$,
we get
$$\s_2 \s_1 \equiv (\bar{B}, b^{-1}) (A, a) \equiv (A \setminus \bar{B}, a) (\bar{B}, b^{-1}) (\bar{B} \setminus A, a^{-1}).$$
Here, we see that $A \setminus \bar{B}= E +b^{-1}$ and $\bar{B}
\setminus  A=\Sig-( C+ B+b)$, so that
$$\s_2 \s_1 \equiv (E +b^{-1}, a) (B, b) (C+ B +b-a^{\pm 1}, a),$$
where $\deg \s_1=\deg (E +b^{-1}, a)> \deg (B, b)$, $\deg (C+ B
+b-a^{\pm 1}, a)$.

\proclaim {Case 6.2} $c \neq b^{-1}$.
\endproclaim

In this case, $c. \Sig > b. \Sig$, since $\deg \s_1$ is determined
by $c$. Apply Lemma ~2.1 ~(ii) to the equalities $|(\bar {A},
a^{-1})^{-1}(u')|=|(\bar {B}, b^{-1})(u')|=|u|$, that is, $|(\bar
{A}, a)(u')|=|(\bar {B}, b^{-1})(u')|=|u|$, to obtain $|(\bar {B}
\setminus \bar {A}, b^{-1})(u')|=|u|$. But since $c \in  \bar {B}
\setminus \bar {A}$ and $c^{-1} \notin  \bar {B} \setminus \bar
{A}$, we have $b.\Sig > c.\Sig$ by Lemma ~2.2, which contradicts
$c. \Sig > b. \Sig$. Hence this case cannot occur.

\proclaim {Case 7} $a^{\pm 1} \in B$, $b \in A$ and $b^{-1} \notin
A$.
\endproclaim

\proclaim {Case 7.1} $c = b$.
\endproclaim

Applying Case 2.2 to $\s_2 \s_1 \equiv (\bar {B}, b^{-1}) (A, a)$,
we get
$$\s_2 \s_1 \equiv (\bar {B}, b^{-1}) (A, a) \equiv (( A \setminus \bar{B}) + \bar{B}, a) (\bar{B}, b^{-1}) (A \cap \bar {B}, a).$$
From the observation that $( A \setminus \bar{ B}) + \bar{ B} = \Sig -( D +b^{-1})$ and $ A \cap \bar { B}= C - b$, it follows that
$$\s_2 \s_1 \equiv ( D+b^{-1}-a^{\pm 1}, a^{-1}) ( B, b) ( C - b, a),$$
where $\deg \s_1=\deg ( D+b^{-1}-a^{\pm 1}, a^{-1})> \deg ( B, b)$, $\deg ( C - b, a)$.

\proclaim {Case 7.2} $c \neq b$.
\endproclaim

As in Case 6.2, $c. \Sig > b. \Sig$. By Lemma 2.1 (ii) applied to
the equalities $|( A, a^{-1})(u')|=|(\bar {  B},
b^{-1})(u')|=|u|$, we get $|(\bar { B} \setminus  A,
b^{-1})(u')|=|u|$. But since $c^{-1} \in  \bar { B} \setminus  A$
and $c \notin  \bar { B} \setminus  A$, we must have $b.\Sig >
c.\Sig$ by Lemma ~2.2, contrary to the fact $c. \Sig > b. \Sig$.
Hence this case cannot happen.

\proclaim {Case 8} $a^{\pm 1} \in  B$ and $b^{\pm 1} \in  A$.
\endproclaim

Apply Lemma 2.1 (ii) to the equalities $|(\bar{ A},
a^{-1})^{-1}(u')|=|(\bar {  B}, b^{-1})(u')|=|u|$, that is,
$|(\bar{ A}, a)(u')|=|(\bar {  B}, b^{-1})(u')|=|u|$, to obtain
$|(\bar { B} \setminus \bar{ A}, b^{-1})(u')|=|u|$; then
$$\s_2 \s_1 \equiv (\bar { B}, b^{-1}) (\bar { A}, a^{-1}) \equiv [(\bar { A} \cap \bar { B}, b^{-1}) (\bar { B} \setminus \bar { A}, b^{-1})] (\bar { A}, a^{-1}).$$
Since $(\bar { B} \setminus \bar { A}, b^{-1}) (\bar { A},
a^{-1})=(\bar { A}, a^{-1})(\bar { B} \setminus \bar { A},
b^{-1})$ by Case 1.1.2, we have
$$\s_2 \s_1 \equiv (\bar { A} \cap \bar { B}, b^{-1}) (\bar { A}, a^{-1}) (\bar { B} \setminus \bar { A}, b^{-1}).$$
It follows from $\bar { A} \cap \bar { B}=\Sig-( C+ B)$ and $\bar { B} \setminus \bar { A}= C - b^{\pm 1}$ that
$$\s_2 \s_1 \equiv ( C +  B - b^{\pm 1}, b)( A, a)( C - b^{\pm 1}, b^{-1}),$$
where $\deg \s_1=\deg ( C +  B - b^{\pm 1}, b)= \deg ( A, a)=\deg ( C - b^{\pm 1}, b^{-1})$.

The proof of the lemma is now completed.
\qed
\enddemo

\proclaim {Remark} The proof of Lemma ~3.1 can be applied without
further change if we replace consideration of a single cyclic word
$u$, the length $|u|$ of $u$, and the total number of occurrences
of $x_i^{\pm 1}$ in $u$ with consideration of a finite sequence
$(u_1, \dots, u_m)$ of cyclic words, the sum $\sum\limits_{i=1}^m
|u_i|$ of the lengths of $u_1, \dots, u_m$, and the total number
of occurrences of $x_i^{\pm 1}$ in $(u_1, \dots, u_m)$,
respectively.
\endproclaim

\proclaim {Lemma 3.2} Let $u$ be a cyclic word in $F_n$ that
satisfies Hypotheses 1.1 and 1.3. Let $\s_i$, $i=1, \dots, \ell$,
be Whitehead automorphisms of the second type such that $|\s_i
\cdots \s_1(u)|=|u|$ for all $i$. Suppose that $\max\limits_{1 \le
i \le \ell} {\deg \s_i} = 1$. Then there exist Whitehead
automorphisms $\t_1, \t_2, \dots, \t_s$ of the second type such
that
$$\s_\ell \cdots \s_2 \s_1(u)= \t_s \cdots \t_2 \t_1(u),$$
where $1 \ge \deg \t_s \ge \deg \t_{s-1} \ge \cdots \ge \deg \t_1$, and $|\t_j \cdots \t_1(u)|=|u|$ for all $j=1, \dots, s$.
\endproclaim

\demo {Proof} We proceed by induction on $\ell$. The case for
$\ell=2$ is already proved in Lemma 3.1. Now let $\s_i$, $i=1,
\dots, \ell+1$, be Whitehead automorphisms of $F_n$ such that
$|\s_i \cdots \s_1(u)|=|u|$ for all $i$ and such that
$\max\limits_{1 \le i \le \ell+1} {\deg \s_i} = 1$. Then by the
induction hypothesis, there exist Whitehead automorphisms $\t_1,
\t_2, \dots, \t_s$ of $F_n$ such that
$$\s_{\ell+1} \s_\ell \cdots \s_2 \s_1(u)= \s_{\ell+1} \t_s \cdots \t_2 \t_1(u),
\tag 3.2
$$
where $1 \ge \deg \t_s \ge \deg \t_{s-1} \ge \cdots \ge \deg \t_1$, and $|\t_j \cdots \t_1(u)|=|u|$ for all $j=1, \dots, s$.

Put $\t_j=( A_j, a_j)$ for $j=1, \dots, s$, and put $\s_{\ell+1}=(
B, b)$. If $\deg \s_{\ell+1}=1$ or $\deg \t_j=0$ for all $j$, then
there is nothing to prove. So let $\deg \s_{\ell+1}=0$, and let
$t$ ($1 \le t \le s$) be such that $\deg \t_s = \deg \t_{s-1}=
\cdots = \deg \t_t=1$ and $\deg \t_{t-1}= \cdots =\deg \t_2 = \deg
\t_1=0$. Upon replacing $\t_i$ and $\s_{\ell+1}$ by $(\bar  A_i,
a_i^{-1})$ and $(\bar  B, b^{-1})$, respectively, if necessary, we
may assume that $x_1 \in  A_i$ for all $t \le i \le s$ and that
$x_1^{\pm 1} \notin  B$. We may also assume without loss of
generality that $(B, b)$ cannot be decomposed into $( B_2, b)(
B_1, b)$, where $B= B_1+ B_2$, $\deg (B_1, b)=\deg (B_2, b)=0$ and
$|(B_1, b) \t_s \cdots \t_1(u)|=|u|$.

\proclaim {Claim 1} We may further assume that $\t_s=(A_s, a_s)$
cannot be decomposed into $(A_{s2}, a_s)( A_{s1}, a_s)$, where $
A_s=A_{s1}+ A_{s2}$, $\deg (A_{s1}, a_s)=0$, $\deg (A_{s2},
a_s)=1$, $|(A_{s1}, a_s) \t_{s-1} \cdots \t_1(u)|=|u|$, and
$a_i^{\pm 1} \notin A_{s1}$ for all $i$ with $t \le i < s$.
\endproclaim

\demo {Proof of Claim 1} Suppose that $\t_s$ can be decomposed in
the same way as in the statement of the claim. Then continuously
applying Case 1 or Case 4 of Lemma ~3.1 to $(A_{s1}, a_s) \t_{s-1}
\cdots \t_t$ at most $1+2+2^2+\cdots+2^{s-t-1}$ times (here, note
that if $s=t$, we do not need to apply Lemma ~3.1), we get
$$( A_{s1}, a_j) \t_{s-1} \cdots \t_t = \t_{s-1}' \cdots \t_t' \e_p \cdots \e_1,
$$
where $\t_{s-1}', \dots, \t_t'$ are Whitehead automorphisms of
degree $1$ and $\e_p, \dots, \e_1$ are Whitehead automorphisms of
degree $0$, so that
$$(B, b) \t_s \cdots \t_t \cdots \t_1(u)= (B, b)( A_{s2}, a_s) \t_{s-1}' \cdots \t_t' \e_p \cdots \e_1 \t_{t-1} \cdots
\t_1(u),
\tag 3.3
$$
where the length of $u$ is constant throughout both chains. We
then replace the chain on the right-hand side of (3.2) with that
of (3.3). \qed
\enddemo

We consider three cases corresponding to whether or not $b=x_1^{\pm 1}$.

\proclaim {Case 1} $b \neq x_1^{\pm 1}$.
\endproclaim
For all $i$ with $t \le i \le s$,  either $b^{\pm 1} \in  A_i$ or
$b^{\pm 1} \notin  A_i$, since $\deg \t_i=1$. If $a_s^{\pm 1} \in
B$, then the required result follows immediately from Case 5 or
Case 8 of Lemma ~3.1 applied to $( B, b) \t_s$. So let $a_s^{\pm
1} \notin  B$. If $b^{\pm 1} \notin  A_s$ and $ A_s \cap  B =
\emptyset$, then by Case 1.1.2 of Lemma 3.1 we have $( B, b) \t_s
\equiv \t_s ( B, b)$. Also if $b^{\pm 1} \in  A_s$ and $ B \subset
A_s$, then Case ~4 of Lemma ~3.1 yields that $( B, b) \t_s \equiv
\t_s ( B, b)$. Hence, in either case, we have
$$( B, b) \t_s \cdots \t_t \cdots \t_1(u)= \t_s ( B, b) \t_{s-1} \cdots \t_t \cdots \t_1(u);
$$
then the desired result follows by induction on $s-t$. Now suppose
that either both $b^{\pm 1} \notin  A_s$ and $ A_s \cap  B \neq
\emptyset$ or both $b^{\pm 1} \in  A_s$ and $ B \nsubseteq  A_s$.
We argue two cases separately.

\proclaim {Case 1.1} $a_s^{\pm 1} \notin  B$, $b^{\pm 1} \notin  A_s$ and $ A_s \cap  B \neq \emptyset$.
\endproclaim
By Case 1.2.3 of Lemma ~3.1, we have $( B, b) \t_s \equiv ( A_s
\setminus  B, a_s) ( B, b) ( A_s \cap  B, a_s)$; thus
$$( B, b) \t_s \cdots \t_t \cdots \t_1(u)= ( A_s \setminus  B, a_s) ( B, b) ( A_s \cap  B, a_s) \t_{s-1} \cdots \t_t \cdots \t_1(u).
$$
By Claim 1, there is $j$ with $t \le j < s$ such that $a_j^{\pm 1}
\in  A_s \cap  B$. Let $r$ be the largest such index.

First suppose that there exists a chain $\eta_m \cdots \eta_1$ of
Whitehead automorphisms $\eta_i=( G_i, g_i)$ of degree 1 with
$g_i^{\pm 1} \notin  B$, $ G_i \subset  A_s$ and $ G_i \cap  B=
\emptyset$ such that $|\eta_i \cdots \eta_1 \t_s \cdots
\t_1(u)|=|u|$ for all $i=1, \dots, m$ and such that $|( H,
a_r^{-1})\eta_m \cdots \eta_1 \t_s \cdots \t_1(u)|=|u|$ for some
Whitehead automorphism $( H, a_r^{-1})$ of degree 1 with $ H
\subset  A_s$. Then
$$\split
(B, b) \t_s \cdots \t_1(u)&=( B, b)\eta_1^{-1} \cdots \eta_m^{-1} \eta_m \cdots \eta_1 \t_s \cdots \t_1(u)\\
&=\eta_1^{-1} \cdots \eta_m^{-1} ( B, b) \eta_m \cdots \eta_1 \t_s
\cdots \t_1(u) \quad \text {by Case 1.1.2 of Lemma 3.1}.
\endsplit
$$
Put $v= \eta_m \cdots \eta_1 \t_s \cdots \t_1(u)$. By Lemma ~2.1
~(ii) applied to $|(\bar B, b^{-1})(v)|=|(H, a_r^{-1})(v)|=|u|$,
we have $|(\bar  B \setminus  H, b^{-1})(v)|=|u|$. It follows from
$\bar B \setminus H = \Sig - (B \cup H)-b^{\pm 1}$ that $|(B \cup
H, b)(v)|=|u|$, so that
$$(B, b) \t_s \cdots \t_1(u)= \eta_1^{-1} \cdots \eta_m^{-1} (H \setminus B, b^{-1}) (B \cup H, b) \eta_m \cdots \eta_1 \t_s \cdots \t_1(u),
$$
where $\deg \eta_i^{-1}=\deg ( H \setminus  B, b^{-1})=\deg (B
\cup  H, b)=\deg \eta_i=1$, as required.

Next suppose that there does not exist such a chain $\eta_m \cdots
\eta_1$ as above. Considering all the assumptions and the
situations above, we can observe that this can possibly happen
only in the case where all of $a_s$ and $a_s^{-1}$ that are lost
in passing from $\t_{s-1} \cdots \t_1(u)$ to $\t_s \cdots \t_1(u)$
were newly introduced in passing from $\t_{q-1} \cdots \t_1(u)$ to
$\t_q \cdots \t_1(u)$ for some $r < q < s$, and where for such
$\t_q=( A_q, a_s^{-1})$ (here note that $a_q=a_s^{-1}$),
$$( B, b) \t_s \cdots \t_t \cdots \t_1(u)= (B, b) (A_s \setminus B, a_s) \t_{s-1} \cdots \t_{q+1} (A_q \setminus (A_s \cap  B), a_s^{-1}) \t_{q-1} \cdots \t_t \cdots \t_1(u),
$$
where the length of $u$ is constant throughout the chain on the
right-hand side. It then follows from Case ~1.1.2 of Lemma ~3.1
applied to $( B, b) ( A_s \setminus  B, a_s)$ that
$$( B, b) \t_s \cdots \t_t \cdots \t_1(u)=
(A_s \setminus B, a_s)(B, b) \t_{s-1} \cdots \t_{q+1} (A_q \setminus (A_s \cap  B), a_s^{-1}) \t_{q-1} \cdots \t_t \cdots \t_1(u).
$$
Then induction on $s-t$ yields the desired result, which completes the proof of Case ~1.1.

\proclaim {Case 1.2} $a_s^{\pm 1} \notin  B$, $b^{\pm 1} \in  A_s$ and $ B \nsubseteq  A_s$.
\endproclaim
In this case, replace $\t_i$ by $(\bar  A_i, a_i^{-1})$ for all $t \le i \le s$ and then follow the arguments of Case 1.1.

\proclaim {Case 2} $b=x_1$.
\endproclaim
We divide this case into two cases according to whether $a_s^{\pm
1} \in B$ or not.

\proclaim {Case 2.1} $a_s^{\pm 1} \in B$.
\endproclaim
In this case, we have by Case 7.1 of Lemma 3.1 applied to $(B,
x_1) \t_s$ that
$$(B, x_1) \t_s \cdots \t_1(u) = (B \setminus  A_s +x_1^{-1}-a_s^{\pm 1}, a_s^{-1}) (B, x_1) (A_s \setminus B -x_1, a_s) \t_{s-1} \cdots \t_1(u).
\tag 3.4
$$
Here if $A_s \setminus B - x_1 =\emptyset$, then
$$(B, x_1) \t_s \cdots \t_t \cdots \t_1(u) = (B - A_s +x_1^{-1}-a_s^{\pm 1}, a_s^{-1}) (B, x_1) \t_{s-1} \cdots \t_t \cdots \t_1(u);$$
hence the desired result follows by induction on $s-t$.

So let $A_s \setminus B - x_1 \neq \emptyset$. By Claim 1, there
is $j$ with $t \le j < s$ such that $a_j^{\pm 1} \in A_s \setminus
B -x_1$. Let $r$ be the largest such index. The following Claims
~2--4 show that we may assume that $a_r$, $a_s$ and $x_1$ belong
to distinct connected components of the dependence graph $\G_u$ of
~$u$.

\proclaim {Claim 2} $a_r$ and $x_1$ belong to distinct connected components of $\G_u$.
\endproclaim

\demo {Proof of Claim 2} Suppose on the contrary that $a_r$ and
$x_1$ belong to the same connected component $C_1$. Put $\Cal W=
\{\a:$ $\a$ is a Whitehead automorphism of degree 0 such that
$|\a(v)|=|v|=|u|$ for some $v \in \text { Orb}_{\text {
Aut}F_n}(u)\}$. Then by (3.4), $( A_s \setminus B-x_1, a_s) \in
\Cal W$ and $( B, x_1) \in \Cal W$. Since $x_1^{\pm 1} \notin A_s
\setminus  B-x_1$ and $a_r^{\pm 1} \in A_s \setminus B-x_1$, we
see from the construction of $\G_u$ that $a_s$ also belongs to
$C_1$ and that every path from $a_r$ or $a_r^{-1}$ to $x_1$ or
$x_1^{-1}$ passes through $a_s$ or $a_s^{-1}$. Also since
$a_r^{\pm 1} \notin  B$ and $a_s^{\pm 1} \in  B$, every path from
$a_s$ or $a_s^{-1}$ to $a_r$ or $a_r^{-1}$ passes through $x_1$ or
$x_1^{-1}$, which contradicts the above fact that every path from
$a_r$ or $a_r^{-1}$ to $x_1$ or $x_1^{-1}$ passes through $a_s$ or
$a_s^{-1}$. \qed
\enddemo

\proclaim {Claim 3} We may assume that $a_s$ and $x_1$ belong to distinct connected components of $\G_u$.
\endproclaim

\demo {Proof of Claim 3} Suppose that $a_s$ and $x_1$ belong to
the same connected component $C_1$. First consider the case where
there exists a chain $\zeta_k \cdots \zeta_1$ of Whitehead
automorphisms $\zeta_i=( E_i, e_i)$ of degree 1 with $e_i^{\pm 1}
\in  B$ and $ E_i \subset ( B + x_1)$ such that $|\zeta_i \cdots
\zeta_1 \t_s \cdots \t_1(u)|=|u|$ for all $i=1, \dots, k$ and such
that $|( H, a_r^{-1}) \zeta_k \cdots \zeta_1 \t_s \cdots
\t_1(u)|=|u|$ for some Whitehead automorphism $( H, a_r^{-1})$ of
degree ~1 with $ H \subset  A_s$. Then
$$\split
(B, x_1) \t_s \cdots \t_1(u)&=(B, x_1)\zeta_1^{-1} \cdots \zeta_k^{-1} \zeta_k \cdots \zeta_1 \t_s \cdots \t_1(u)\\
&=\rho_k \cdots \rho_1 (B, x_1) \zeta_k \cdots \zeta_1 \t_s \cdots
\t_1(u) \quad \text {by Case 7.1 of Lemma 3.1},
\endsplit
$$
where $\rho_i=(B \setminus E_{k+1-i} +x_1^{-1}-e_{k+1-i}^{\pm 1},
e_{k+1-i}^{-1})$ for $i=1, \dots, k$. Put $v= \zeta_k \cdots
\zeta_1 \t_s \cdots \t_1(u)$. Then $|(B, x_1)(v)|=|(H,
a_r^{-1})(v)|=|u|$, that is, $|(B, x_1)(v)|=|(\bar H,
a_r)(v)|=|u|$. By Lemma ~2.1 ~(ii) applied to these equalities, we
have $|(\bar H \setminus B, a_r)(v)|=|u|$, so that
$$|(H + (\bar H \setminus B), a_r) ( H, a_r^{-1}) \zeta_k \cdots \zeta_1 \t_s \cdots \t_1(u)|=|u|.$$
It then follows from $H+(\bar H \setminus B) = \Sig-(B \setminus
H)-a_r^{\pm 1}$ that
$$|(B \setminus H, a_r^{-1}) (H, a_r^{-1})\zeta_k \cdots \zeta_1 \t_s \cdots \t_1(u)|=|u|.$$
This implies that $(B \setminus H, a_r^{-1}) \in \Cal W$, where
$\Cal W$ is defined in the proof of Claim ~2. Since $a_s^{\pm 1}
\in  B \setminus  H$ and $x_1^{\pm 1} \notin  B \setminus  H$,
$a_r$ must also belong to $C_1$ by the construction of $\G_u$,
which contradicts Claim ~2.

Next consider the case where there does not exist such a chain
$\zeta_k \cdots \zeta_1$ as above. Considering all the assumptions
and the situations above, we can observe that this can possibly
happen only in the case where all of $a_s$ and $a_s^{-1}$ that are
lost in passing from $\t_{s-1} \cdots \t_1(u)$ to $\t_s \cdots
\t_1(u)$ were newly introduced in passing from $\t_{q-1} \cdots
\t_1(u)$ to $\t_q \cdots \t_1(u)$ for some $r < q < s$, and where
for such $\t_q=( A_q, a_s^{-1})$ (here note that $a_q=a_s^{-1}$),
$$(B, x_1) \t_s \cdots \t_t \cdots \t_1(u)= (B, x_1)(A_s \cap B, a_s) \t_{s-1} \cdots \t_{q+1} (A_q \setminus (A_s \setminus B), a_s^{-1}) \t_{q-1} \cdots \t_t \cdots \t_1(u),
$$
where the length of $u$ is constant throughout the chain on the
right-hand side. It then follows from Case ~7.1 of Lemma ~3.1
applied to $( B, x_1) ( A_s \cap  B, a_s)$ that
$$\multline
(B, x_1) \t_s \cdots \t_t \cdots \t_1(u) \\
= (B \setminus  A_s +x_1^{-1} -a_s^{\pm 1}, a_s^{-1})( B, x_1)
\t_{s-1} \cdots \t_{q+1} ( A_q \setminus ( A_s \setminus  B),
a_s^{-1}) \t_{q-1} \cdots \t_t \cdots \t_1(u).
\endmultline
$$
So in this case, apply induction on $s-t$ to get the desired result of the lemma, which completes the proof of Claim ~3.
\qed
\enddemo

\proclaim {Claim 4} $a_r$ and $a_s$ belong to distinct connected components of $\G_u$.
\endproclaim

\demo {Proof of Claim 4} Suppose on the contrary that $a_r$ and
$a_s$ belong to the same connected component. Note that $a_r^{\pm
1} \notin B$, $a_s^{\pm 1} \in B$ and that $(B, x_1) \in \Cal W$,
where $\Cal W$ is defined in the proof of Claim ~2. It then
follows from the construction of $\G_u$ that $a_s$ and $x_1$ must
belong to the same connected component, which contradicts Claim
~3. \qed
\enddemo

So let $C_1$, $C_{r'}$ and $C_{s'}$ be the distinct connected
components of $\G_u$ containing $x_1$, $a_r$, and $a_s$ in that
order. Here notice that $C_1$ consists of only $x_1^{\pm 1}$,
since there exists a Whitehead automorphism $(A_s, a_s)$ of degree
1 such that $a_s \notin C_1$ and such that $|(A_s,
a_s)(v)|=|v|=|u|$ for some $v \in \text {Orb}_{\text {Aut}
F_n}(u)$ (see Remark ~(iii) in the Introduction).

Put $u_1=\t_{t-1} \cdots \t_1(u)$.

\proclaim{Claim 5} We may assume that $\t_i\t_j \equiv \t_j\t_i$
for all $1 \le i \neq j \le t-1$.
\endproclaim

\demo{Proof of Claim 5} Put $\Cal M=\{v : v=\phi (u) \ \text{and}
\ |v|_{C_i}=|u|_{C_i} \ \text{for all} \ i=1, \dots, n$, where
$\phi$ is a chain of Whitehead automorphisms of degree 0
throughout which the length of $u$ is constant$\}$. Taking an
appropriate $v \in \Cal M$, we have Whitehead automorphisms
$\d_j=(D_j, d_j)$ of $F_n$ of degree 0 such that
$$u_1=\d_h \cdots \d_1(v),
\tag 3.5
$$
where $|\d_j \cdots \d_1(v)|=|v|$ and $|\d_j \cdots
\d_1(v)|_{C_{k_j}}
> |v|_{C_{k_j}}$ for the connected component $C_{k_j}$ containing $d_j$
and for each $j=1, \dots, h$. Then for any $\d_i=(D_i, d_i)$ and
$\d_j=(D_j, d_j)$ with $d_j \neq d_i^{\pm 1}$, if we replace
$\d_i$ and $\d_j$ with $(\bar D_i, d_i^{-1})$ and $(\bar D_j,
d_j^{-1})$, respectively, if necessary so that $d_i^{\pm 1} \notin
D_j$ and $d_j^{\pm 1} \notin D_i$, then $D_i \cap  D_j =
\emptyset$. Hence by Case ~1.1.2 of Lemma ~3.1 that $\d_j \d_i
\equiv \d_i \d_j$; thus (3.5) can be re-written as
$$u_1=\d_{p t_p}^{q_{p t_p}} \cdots \d_{p1}^{q_{p 1}} \cdots \d_{1t_1}^{q_{1 t_1}} \cdots \d_{11}^{q_{11}}(v),
\tag 3.6
$$
where $d_{k i}=d_{k i'}$ and $D_{k i} \neq  D_{k i'}$ provided $i
\neq i'$; $d_{k' i} \neq d_{k i}^{\pm 1}$ and $(\d_{k'
t_{k'}}^{q_{k' t_{k'}}} \cdots \d_{k'1}^{q_{k'1}})(\d_{k
t_k}^{q_{k t_k}} \cdots \d_{k 1}^{q_{k 1}}) \equiv (\d_{k
t_k}^{q_{k t_k}} \cdots \d_{k 1}^{q_{k 1}}) (\d_{k' t_{k'}}^{q_{k'
t_{k'}}} \cdots \d_{k'1}^{q_{k'1}})$ provided $k \neq k'$. Here we
may assume by Case ~1.2.1 of Lemma ~3.1 that $D_{k i} \subset D_{k
i'}$ if $i < i'$. Then $\d_{k i'}\d_{k i} \equiv \d_{k i}\d_{k
i'}$ by Case ~1.2.1 of Lemma ~3.1; hence $\d_{k' i'} \d_{ki}
\equiv \d_{ki}\d_{k'i'}$ for any $\d_{ki}$ and $\d_{k'i'}$ in
chain (3.6). Thus replace $\t_{t-1} \cdots \t_1(u)$ with the
right-hand side of (3.6) to get our desired result. \qed
\enddemo

By Claim 5, we may write
$$u_1=\t_{t-1} \cdots \t_p \t_{p-1} \cdots \t_1(u),$$
where $\t_i$ has multiplier in $C_{r'}$ provided $p \le i \le
t-1$; $\t_i$ has multiplier not in $C_{r'}$ provided $1 \le i \le
p-1$. Put
$$u_2=\t_{p-1} \cdots \t_1(u).$$
Note that the the number of $C_{r'}$-syllables of $u$ remains
unchanged throughout this chain.

\proclaim {Claim 6} There exist Whitehead automorphisms
$\e_i=(E_i, a_i)$, $t \le i \le s$, such that $|\e_i \cdots \e_t
(u_2)|=|u|$ for all $i=t, \dots, s$, where $E_i=\emptyset$
provided $a_i \in C_{r'}$; $E_i$ is one of the three forms $A_i$,
$A_i + C_{r'}$ and $A_i - C_{r'}$, whichever is smallest possible
with priority given to lower $i$, provided $a_i \notin C_{r'}$.
\endproclaim

\demo {Proof of Claim 6} Suppose the contrary. It can possibly
happen only when the number of $C_{r'}$-syllables of $u_2$ is
decreased by $\t_j \cdots \t_t \t_{t-1} \cdots \t_p$ (for some $j
\ge t$) followed by a chain of Whitehead automorphisms of degree 0
with multiplier in $C_{r'}$, where the length of $u_2$ is constant
throughout the chain. Choosing the smallest such index $j$, put
$\{j_1, \dots, j_k\}=\{i: t \le i \le j$ and $\t_i$ has multiplier
in $C_{r'}\}$. Then we can observe that there is a chain $\zeta_m
\cdots \zeta_1$ of Whitehead automorphisms of degree ~0 with
multiplier in $C_{r'}$ such that $|\zeta_m \cdots \zeta_1 \t_{j_k}
\cdots \t_{j_1} \t_{t-1} \cdots \t_p(u_2)|=|u_2|$ and the number
of $C_{r'}$-syllables of $\zeta_m \cdots \zeta_1 \t_{j_k} \cdots
\t_{j_1} \t_{t-1} \cdots \t_p(u_2)$ is less than that of $u_2$.
This is a contradiction, because through the chain $\zeta_m \cdots
\zeta_1 \t_{j_k} \cdots \t_{j_1} \t_{t-1} \cdots \t_p$ only
$C_1$-syllables and $C_{r'}$-syllables can mix and increasing the
number of $C_1$-syllables cannot reduce the number of
$C_{r'}$-syllables. \qed
\enddemo

For the chain $\e_s \cdots \e_t$, we consider two cases
separately.

\proclaim {Case 2.1.1} $|(B, x_1) \e_s \cdots \e_t (u_2)|=|u|$.
\endproclaim
For the Whitehead automorphisms $\d_i=(D_i, d_i)$ ($p \le i < t$),
where $D_i=A_i \setminus B$ and $d_i=a_i$ provided $x_1^{\pm 1}
\notin A_i$; $D_i=\bar{A_i} \setminus B$ and $d_i=a_i^{-1}$
provided $x_1^{\pm 1} \in A_i$, and for the Whitehead
automorphisms $\o_j=(F_j, a_{t+s-j}^{-1})$ and $\nu_j=(H_j, a_j)$
($t \le j \le s$), where $F_j = \emptyset$ provided $a_{t+s-j} \in
C_{r'}+B$; $F_j=E_{t+s-j} \setminus B$ provided $a_{t+s-j} \notin
C_{r'}+B$; $H_j=\emptyset$ provided $a_j \in B$; $H_j=A_j
\setminus B$ provided $a_j \notin B$, we have
$$
(B, x_1) \t_s \cdots \t_1(u)= \nu_s \cdots \nu_t \d_{t-1} \cdots
\d_p \o_s \cdots \o_t (B, x_1) \e_s \cdots \e_t \t_{p-1} \cdots
\t_1(u), \tag 3.7
$$
where the length of $u$ is constant throughout the chain on the
right-hand side. By Case ~1, it suffices to consider only the
chain $(B, x_1) \e_s \cdots \e_t \t_{p-1} \cdots \t_1(u)$. Since
for every $j$ either $\deg \e_j=1$ or $\e_j=1$ and since $\e_r=1$,
the desired result follows by induction on $s-t$ from (3.7).

\proclaim {Case 2.1.2} $|(B, x_1) \e_s \cdots \e_t (u_2)|>|u|$.
\endproclaim
We see that this case can possibly happen only when the cyclic
word $\e_s \cdots \e_t (u_2)$ contains a subword of the form $(x_1
w_1 w_2 w_3)^{\theta}$, where $\theta={\pm 1}$, $w_1$ ($w_1$ may
be the empty word), $w_2$ and $w_3$ are words in $B$, $C_{r'}$ and
$C_{s'}$, respectively, and not all of the letters in $w_3$ were
newly introduced in passing from $u_2$ to $\e_s \cdots \e_t(u_2)$.

By Claim~5, we may write
$$u_1=\t_{t-1} \cdots \t_q \t_{q-1} \cdots \t_1(u),$$
where $\t_i$ has multiplier in $C_{s'}$ provided $q \le i \le
t-1$; $\t_i$ has multiplier not in $C_{s'}$ provided $1 \le i \le
q-1$. Put
$$u_3=\t_{q-1} \cdots \t_1(u).$$
Notice that the the number of $C_{s'}$-syllables of $u$ remains
unchanged throughout this chain.

\proclaim {Claim 7} There exist Whitehead automorphisms
$\lambda_i=(J_i, a_i)$, $t \le i \le s$, such that $|\lambda_i
\cdots \lambda_t (u_3)|=|u|$ for all $i=t, \dots, s$, where $
J_i=\emptyset$ provided $a_i \in C_{s'}$; $ J_i$ is one of the
three forms $A_i$, $A_i + C_{s'}$ and $ A_i - C_{s'}$, whichever
is largest possible with priority given to lower $i$, provided
$a_i \notin C_{s'}$.
\endproclaim

\demo {Proof of Claim 7} Suppose the contrary. In view of all the
assumptions and the situations above, this can possibly happen
only when the number of $C_{s'}$-syllables of $u_3$ is decreased
by $\t_j \cdots \t_t \t_{t-1} \cdots \t_q$ (for some $j \ge t$)
followed by a chain of Whitehead automorphisms of degree 0 with
multiplier in $C_{s'}$, where the length of $u_3$ is constant
throughout the chain. Choosing the smallest such index $j$, put
$\{j_1, \dots, j_k\}=\{i: t \le i \le j$ and $\t_i$ has multiplier
in $C_{s'}\}$. Then we can observe that there exists a chain $\d_m
\cdots \d_1$ of Whitehead automorphisms of degree ~0 with
multiplier in $C_{s'}$ such that $|\d_m \cdots \d_1 \t_{j_k}
\cdots \t_{j_1}\t_{t-1} \cdots \t_q (u_3)|=|u|$, and such that the
number of $C_{s'}$-syllables of $\d_m \cdots \d_1 \t_{j_k} \cdots
\t_{j_1}\t_{t-1} \cdots \t_q (u_3)$ is less than that of $u_3$.
Reasoning as in Claim ~6, we get a contradiction, which completes
the proof of Claim ~7. \qed
\enddemo

We then see that $|(B, x_1) \lambda_s \cdots \lambda_t
(u_3)|=|u|$. Furthermore, for the Whitehead automorphisms
$\d_i=(D_i, d_i)$ ($q \le i < t$), where $D_i=A_i \cap B$ and
$d_i=a_i$ provided $x_1^{\pm 1} \notin A_i$; $D_i=\bar{A_i} \cap
B$ and $d_i=a_i^{-1}$ provided $x_1^{\pm 1} \in A_i$, and for the
Whitehead automorphisms $\o_j=(K_j, a_{t+s-j})$ and $\nu_j=(H_j,
a_j^{-1})$ ($t \le j \le s$), where $K_j = \emptyset$ provided
$a_{t+s-j} \notin B-C_{s'}$; $K_j= B \setminus J_{t+s-j} +
x_1^{-1}-a_{t+s-j}^{\pm 1}$ provided $a_{t+s-j} \in B-C_{s'}$;
$H_j=\emptyset$ provided $a_j \notin B$; $H_j= B \setminus
A_j+x_1^{-1}-a_j^{\pm 1}$ provided $a_j \in B$,
$$
(B, x_1) \t_s \cdots \t_1(u)= \nu_s \cdots \nu_t \d_{t-1} \cdots
\d_q \o_s \cdots \o_t (B, x_1) \lambda_s \cdots \lambda_t \t_{q-1}
\cdots \t_1(u), \tag 3.8
$$
where the length of $u$ is constant throughout the chain on the
right-hand side. By Case ~1, it suffices to consider only the
chain $(B, x_1) \lambda_s \cdots \lambda_t \t_{q-1} \cdots
\t_1(u)$. Since for every $j$ either $\deg \lambda_i=1$ or
$\lambda_i=1$ and since $\lambda_s=1$, the desired result follows
by induction on $s-t$ from (3.8). This completes the proof of Case
~2.1.2.

\proclaim {Case 2.2} $a_s^{\pm 1} \notin  B$.
\endproclaim
In this case, replace $(B, x_1)$ and $\t_i$ by $(\bar  B,
x_1^{-1})$ and $(\bar A_i, a_i^{-1})$ for all $t \le i \le s$,
respectively, and then follow the arguments of Case 2.1.

\proclaim {Case 3} $b=x_1^{-1}$.
\endproclaim

Replace $( B, x_1^{-1})$ by $(\bar  B, x_1)$ and then repeat the arguments of Case ~2.
\qed
\enddemo

\proclaim {Remark} The proof of Lemma ~3.2 can be applied without
further change if we replace consideration of a single cyclic word
$u$, the length $|u|$ of $u$, the total number of occurrences of
$x_j^{\pm 1}$ in $u$, and the $C_j$-syllable length $|u|_{C_j}$
with consideration of a finite sequence $(u_1, \dots, u_m)$ of
cyclic words, the sum $\sum\limits_{i=1}^m |u_i|$ of the lengths
of $u_1, \dots, u_m$, the total number of occurrences of $x_j^{\pm
1}$ in $(u_1, \dots, u_m)$, and the sum $\sum\limits_{i=1}^m
|u_i|_{C_j}$ of the $C_j$-syllable lengths of $u_1, \dots, u_m$,
respectively.
\endproclaim

We are now in a position to prove Theorem ~1.4.

\proclaim{Proof of Theorem 1.4}
\endproclaim

The proof proceeds by double induction on $\ell$ and $r$, where
$\ell$ is the length of the chain $\s_{\ell} \cdots \s_2 \s_1$ and
$r=\max\limits_{1 \le i \le \ell} \deg \s_i$. The base steps were
already proved in Lemma ~3.1 (the case for $\ell=2$ and any $r$)
and Lemma ~3.2 (the case for $r=1$ and any $\ell$).

Let $\s_i$, $i=1, \dots, \ell+1$ ($\ell+1 \ge 3$), be Whitehead
automorphisms of $F_n$ such that $|\s_i \cdots \s_1(u)|=|u|$ for
all $i=1, \dots, \ell+1$ and such that ${\max\limits_{1 \le i \le
\ell+1} \deg \s_i}=r+1 \ge 2$. By the induction hypothesis on
$\ell$, there exist Whitehead automorphisms $\t_1, \t_2, \dots,
\t_s$ of $F_n$ such that
$$\s_{\ell+1} \s_\ell \cdots \s_2 \s_1(u)= \s_{\ell+1} \t_s \cdots \t_2 \t_1(u),
$$
where $r+1 \ge \deg \t_s \ge \deg \t_{s-1} \ge \cdots \ge \deg
\t_1$, and $|\t_j \cdots \t_1(u)|=|u|$ for all $j=1, \dots, s$.

If either $\deg \s_{\ell+1}=r+1$ or both $\deg \t_s \le r$ and
$\deg \s_{\ell+1} \ge r$, then there is nothing to prove. Also if
$\deg \t_s \le r$ and $\deg \s_{\ell+1}< r$, then we are done by
the induction hypothesis on ~$r$. So let $t$ ($1 \le t \le s$) be
such that $\deg \t_i =r+1$ provided $t \le i \le s$ and $\deg \t_i
\le r$ provided $1 \le i < t$, and let $\deg \s_{\ell+1} \le r$.

Put $\t_j=( A_j, a_j)$ for $j=1, \dots, s$ and $\s_{\ell +1}=( B,
b)$. Upon replacing $\t_i$ and $\s_{\ell +1}$ by $(\bar  A_i,
a_i^{-1})$ and $(\bar  B, b^{-1})$, respectively, if necessary, we
may assume that $x_{r+1} \in  A_i$ for all $t \le i \le s$ and
that $x_{r+1}^{\pm 1} \notin  B$. We may also assume without loss
of generality that $( B, b)$ cannot be decomposed to $( B_2, b)(
B_1, b)$, where $ B= B_1+ B_2$, $\deg (B_1, b)=\deg ( B_2, b)=0$
and $|( B_1, b) \t_s \cdots \t_1(u)|=|u|$. We may further assume
as in Claim ~1 of Lemma ~3.2 that $\t_s=(A_s, a_s)$ cannot be
decomposed to $(A_{s2}, a_s)( A_{s1}, a_s)$, where $A_s= A_{s1}+
A_{s2}$, $\deg (A_{s1}, a_s) \le r$, $\deg (A_{s2}, a_s)=r+1$, $|(
A_{s1}, a_s) \t_{s-1} \cdots \t_1(u)|=|u|$, and $a_i^{\pm 1}
\notin A_{s1}$ for all $i$ with $t \le i < s$.

There are three cases to consider.

\proclaim {Case 1} $b=x_1$.
\endproclaim
If ${a_i}^{\pm 1} \notin  B$ for all $t \le i \le s$, then
continuous application of Cases 1--4 of Lemma ~3.1 to $(B, x_1)
\t_s \cdots \t_t$ at most $1+2+2^2+\cdots+2^{s-t}$ times together
with the induction hypothesis on $r$ yields the desired result.
The following Claim shows that it is indeed true that ${a_i}^{\pm
1} \notin  B$ for all $t \le i \le s$.

\proclaim {Claim} ${a_i}^{\pm 1} \notin  B$ for all $t \le i \le
s$.
\endproclaim

\demo {Proof of the Claim} Suppose on the contrary that
${a_i}^{\pm 1} \in  B$ for some $t \le i \le s$. First let
$a_s^{\pm 1} \in  B$. If either $x_1 \in  A_s$ or $x_1^{-1} \in
A_s$ but not both, then we have a contradiction by Cases 6.2 and
7.2 of Lemma ~3.1, since $\deg \t_s =r+1 \ge 2$. If $x_1^{\pm 1}
\in  A_s$, then by Case ~8 of Lemma ~3.1,
$$( B, x_1)( A_s, a_s) \equiv ( A_s \cup  B -x_1^{\pm 1}, x_1)( A_s, a_s)( A_s \setminus  B-x_1^{\pm 1}, x_1^{-1}),$$
but the existence of $( A_s \setminus  B-x_1^{\pm 1}, x_1^{-1})$
in this chain contradicts Lemma ~2.2, because $x_{r+1} \in  A_s
\setminus  B-x_1^{\pm 1}$ and $x_{r+1}^{-1} \notin  A_s \setminus
B-x_1^{\pm 1}$. If $x_1^{\pm 1} \notin  A_s$, then by Case 5 of
Lemma ~3.1,
$$(B, x_1)(A_s, a_s) \equiv (A_s \setminus  B, x_1^{-1})( A_s, a_s)( A_s \cup  B, x_1),$$
but the existence of $(A_s \cup B, x_1)$ in this chain also
contradicts Lemma ~2.2, since $x_{r+1} \in A_s \cup B$ and
$x_{r+1}^{-1} \notin A_s \cup B$.

Next let $a_s^{\pm 1} \notin  B$. Suppose that ${a_i}^{\pm 1} \in
B$ for some $t \le i < s$. Let $k$ be the largest such index. Put
$v=\t_{k-1} \cdots \t_1(u)$. If $x_1 \in A_k$ and $x_1^{-1} \notin
A_k$, then we can observe based on all the assumptions and the
situations above that there exists a Whitehead automorphism $( F,
x_1)$ of degree $r+1$ with $( B \cup A_k -x_1) \subseteq  F$ such
that $|( F, x_1)\t_k(v)|=|u|$. But this yields a contradiction to
Lemma ~2.2, since $x_{r+1} \in  F$ and $x_{r+1}^{-1} \notin  F$.
For a similar reason, the case where $x_1 \notin  A_k$ and
$x_1^{-1} \in  A_k$ cannot happen, either. So $ A_k$ must contain
either both of $x_1^{\pm 1}$ or none of $x_1^{\pm 1}$.

If there exists a chain $\zeta_p \cdots \zeta_1$ of Whitehead
automorphisms of degree less than or equal to $r+1$ such that $|(
B, x_1) \t_k \zeta_p \cdots \zeta_1(v)|=|\t_k \zeta_p \cdots
\zeta_1(v)|=|\zeta_p \cdots \zeta_1(v)|=|u|$, then as in the case
where $a_s^{\pm 1} \in  B$ we reach a contradiction. Otherwise,
choose chains $\zeta_p \cdots \zeta_1$ and $\o_q \cdots \o_1$ of
Whitehead automorphisms of degree less than or equal to $r+1$ with
$q$ smallest possible such that $|\o_j \cdots \o_1 \t_k \zeta_p
\cdots \zeta_1(v)|=|\t_k \zeta_p \cdots \zeta_1(v)|=|\zeta_p
\cdots \zeta_1(v)|=|u|$ for all $j=1, \dots, q$, and such that $|(
B, x_1) \o_q \cdots \o_1 \t_k \zeta_p \cdots \zeta_1(v)|=|u|$.
Clearly $q \le s-k$.

Put $\o_j=( G_j, g_j)$ for $j=1, \dots, q$. If $x_1^{\pm 1} \notin
A_k$, then we see from the choice of $k$ and the chain $\o_q
\cdots \o_1$ that $g_1^{\pm 1} \notin  A_k$. We also see that for
the Whitehead automorphisms $\g_j=( H_j, g_j)$, $j=1, \dots, q$,
where $ H_j= G_j \setminus  A_k$ provided $a_k^{\pm 1} \notin
G_j$; $ H_j=  G_j \cup  A_k$ provided $a_k^{\pm 1} \in  G_j$, $|(
B, x_1) \g_q \cdots \g_1 \t_k \zeta_p \cdots \zeta_1 (v)|=|\g_j
\cdots \g_1 \t_k \zeta_p \cdots \zeta_1 (v)|=|u|$ for all $j=1,
\dots, q$. Then by Case ~1.1.2 or Case ~5 of Lemma ~3.1, we have
$\g_1 \t_k  \equiv \t_k \g_1$, which means the chain $\g_q \cdots
\g_2$ of shorter length has the same property as $\o_q \cdots
\o_1$ does, contrary to the choice of the chain $\o_q \cdots
\o_1$. If $x_1^{\pm 1} \in  A_k$, replace $\t_k$ by $(\bar  A_k,
a_k^{-1})$. Then we get a contradiction in the same way, which
completes the proof of the claim. \qed
\enddemo

\proclaim {Case 2} $b=x_1^{-1}$.
\endproclaim
Repeat similar arguments to those in Case 1.

\proclaim {Case 3} $b \neq x_1^{\pm 1}$.
\endproclaim
Let $p$ ($1 \le p \le t$) be such that $\deg \t_i =0$ provided $1
\le i < p$; $\deg \t_i \ge 1$ provided $p \le i \le s$. As in
Claim ~5 of Lemma ~3.2, we may assume that $\t_i\t_j \equiv
\t_j\t_i$ for all $1 \le i \neq j < p$. So there exists $q$ with
$1 \le q \le p$ such that $\t_i$ has multiplier in $C_1$ provided
$1 \le i <q$; $\t_i$ has multiplier not in $C_1$ provided $q \le i
<p$.

Put $w=\t_{q-1} \cdots \t_1(u)$. Notice that $C_i$-syllables
remain unchanged throughout the chain $\t_{q-1} \cdots \t_1$ for
all $i \ge 2$. Write
$$w=y_1u_1y_2u_2 \cdots y_mu_m \hskip 0.2in \text {without cancellation},
\tag 3.9
$$
where for each $i=1, \dots, m$, $y_i= x_1$ or $y_i=x_1^{-1}$, and
$u_i$ is a (non-cyclic) subword in $\{x_2, \dots, x_n\}^{\pm 1}$.
Let $F_{n+3}$ be the free group on the set
$$\{x_1, \dots, x_n, x_{n+1}, x_{2n+1}, x_{3n+1}\}.$$ From (3.9)
we construct a sequence $V_{w}=(v_1, v_2, \dots, v_m)$ of cyclic
words $v_1, v_2, \dots, v_m$ in $F_{n+3}$ with $\sum
\limits_{j=1}^m |v_j|=2|u|$, where $m$ is the total number of
occurrences of $x_1^{\pm 1}$ in $u$, as follows: for each $j=1,
\dots, m$,
$$\align
&\text {if $y_j=x_1$ and $y_{j+1}=x_1$, then $v_j=x_1 u_j x_{3n+1} u_j^{-1}$;}\\
&\text {if $y_j=x_1^{-1}$ and $y_{j+1}=x_1$, then $v_j=x_{n+1} u_j x_{3n+1} u_j^{-1}$;} \\
&\text {if $y_j=x_1$ and $y_{j+1}=x_1^{-1}$, then $v_j=x_1 u_j x_{2n+1} u_j^{-1}$;} \\
&\text {if $y_j=x_1^{-1}$ and $y_{j+1}=x_1^{-1}$, then
$v_j=x_{n+1} u_j x_{2n+1} u_j^{-1}$,}
\endalign
$$
where $y_{m+1}=y_1$.

Put $I=\{x_1, x_{n+1}, x_{2n+1}, x_{3n+1}\}^{\pm 1}$. From now on,
when we say that $(S, s)$ is a Whitehead automorphism of
$F_{n+3}$, the following restrictions are imposed on $S$ and $s$:

\roster \item"(1)" $s \in \{x_2, \dots, x_n\}^{\pm 1}$.

\item"(2)" $S$ satisfies one of (i) $I \subseteq S$; (ii) $I \cap
S=\{x_1, x_{2n+1}\}^{\pm 1}$; (iii) $I \cap S=\{x_{n+1},
x_{3n+1}\}^{\pm 1}$; (iv) $I \cap S=\emptyset$.
\endroster
Then we can prove the following

\proclaim{Claim 1} For each Whitehead automorphism $\t=(A, a)$ of
$F_n$ such that $a \neq x_1^{\pm 1}$ and $|\t(w)|=|w|$, there
exists a Whitehead automorphism $\a$ of $F_{n+3}$ such that
$\sum\limits_{j=1}^m |\a(v_j)|=\sum\limits_{j=1}^m |v_j|$ and
$\a(V_{w})=V_{\t(w)}$.
\endproclaim

\demo{Proof of Claim 1} Given a Whitehead automorphism $\t=(A,
a)$, we define a Whitehead automorphism $\a$ of $F_{n+3}$ as
follows: If $x_1^{\pm 1} \in A$, then $\a=(A +x_{n+1}^{\pm
1}+x_{2n+1}^{\pm 1}+x_{3n+1}^{\pm 1}, a)$; if only $x_1 \in A$,
then $\a=(A +x_1^{-1}+x_{2n+1}^{\pm 1}, a)$; if only $x_1^{-1} \in
A$, then $\a=(A-x_1^{-1}+x_{n+1}^{\pm 1}+x_{3n+1}^{\pm 1}, a)$; if
$x_1^{\pm 1} \notin A$, then $\a=(A, a)$.

Then each newly introduced letter $x_r^{\pm 1}$ in passing from
$w$ to $\t(w)$ that remains in $\t(w)$ produces two newly
introduced letters $x_r^{\pm 1}$ in passing from $V_{w}$ to
$\a(V_{w})$ that remain in $\a(V_{w})$, and vice versa. Also each
letter $x_r^{\pm 1}$ in $w$ that is lost in passing from $w$ to
$\t(w)$ produces two letters $x_r^{\pm 1}$ in $V_{w}$ that are
lost in passing from $V_{w}$ to $\a(V_{w})$, and vice versa. This
yields that $\sum\limits_{j=1}^m |\a(v_j)|=\sum\limits_{j=1}^m
|v_j|$.

Moreover it is clear that $\a(V_{w})=V_{\t(w)}$. \qed
\enddemo

The following claim is a converse of Claim ~1.

\proclaim{Claim 2} For each Whitehead automorphism $\a=(S, s)$ of
$F_{n+3}$ such that $\sum\limits_{j=1}^m
|\a(v_j)|=\sum\limits_{j=1}^m |v_j|$, there exists a Whitehead
automorphism $\t=(A, a)$ of $F_n$ such that $a \neq x_1^{\pm 1}$,
$|\t(w)|=|w|$ and such that $\a(V_{w})=V_{\t(w)}$.
\endproclaim

\demo{Proof of Claim 2} Given a Whitehead automorphism $\a=(S,s)$
of $F_{n+3}$, put $T=S \setminus I$. And define a Whitehead
automorphism $\t$ of $F_n$ as follows: $\t=(T+x_1^{\pm 1}, s)$
provided $I \subseteq S$; $\t=(T+x_1, s)$ provided $I \cap
S=\{x_1, x_{2n+1}\}^{\pm 1}$; $\t=(T+x_1^{-1}, s)$ provided $I
\cap S=\{x_{n+1}, x_{3n+1}\}^{\pm 1}$; $\t=(T, s)$ provided $I
\cap S=\emptyset$. Then reasoning in the same way as in Claim ~1,
we get a desired result. \qed
\enddemo

For each $\t_i=(A_i, a_i)$, $q \le i \le s$, define a Whitehead
automorphism $\a_i$ of $F_{n+3}$ as in Claim ~1. Also as in Claim
~1, define a Whitehead automorphism $\b$ of $F_{n+3}$ from
$\s_{\ell+1}=(B, b)$. Then we have $\sum\limits_{j=1}^m |\b \a_s
\cdots \a_q(v_j)|=\sum\limits_{j=1}^m |\a_i \cdots
\a_q(v_j)|=\sum\limits_{j=1}^m |v_j|$ for all $i=q, \dots, s$.
Moreover, by the construction of $\a_i$ and $\b$, the Whitehead
automorphisms $\a_i$ and $\b$ of $F_{n+3}$ are of degree at most
$r+1$, and each of defining sets of $\a_i$ and $\b$ contains
either both of $x_1^{\pm 1}$ or none of $x_1^{\pm 1}$. This yields
the same situation as for a chain of Whitehead automorphisms of
$F_{n+3}$ of maximum degree ~$r$.

Here we notice from Claims ~1 and 2 that if $\G_u$ consists of $g$
connected components, then either $\G_{V_{w}}$ consists of $g+1$
connected components such that $C_i$ equals $C_i$ of $\G_u$ for
all $C_i$'s of $\G_{V_{w}}$ with $C_i \neq C_1$ and $C_i \neq
C_{n+1}$, $C_1$ equals $C_1$ of $\G_u$ plus $x_{2n+1}^{\pm 1}$,
and such that $C_{n+1}=\{x_{n+1}, x_{3n+1}\}^{\pm 1}$; or
$\G_{V_{w}}$ consists of $g$ connected components such that $C_i$
equals $C_i$ of $\G_u$ for all $C_i$'s of $\G_{V_{w}}$ with $C_i
\neq C_1$ and such that $C_1$ equals $C_1$ of $\G_u$ plus
$\{x_{n+1}, x_{2n+1}, x_{3n+1}\}^{\pm 1}$.

The sequence $V_{w}=(v_1, \dots, v_m)$ satisfies neither
Hypothesis ~1.1 nor Hypothesis 1.3. However, this fact does not
affect the proof of the base steps of the induction (that is,
Lemmas 3.1 and 3.2) for the following four reasons: first each of
the Whitehead automorphisms $\a_i$ and $\b$ has multiplier only in
$\{x_2, \dots, x_n\}^{\pm 1}$; second only the proof of Case ~2.1
of Lemma ~3.2 is concerned with the $C_i$-syllable length, but in
the proof of Case ~2.1 $a_r$ or $a_s$ cannot belong to the
connected component $C_1$ of $\G_{V_{w}}$ (in fact, if $a_r$ or
$a_s$ belonged to $C_1$, such a situation as Case ~2.1 could not
occur); third Claim ~5 holds for $V_{w}$ by replacing $\Cal M$
with the set $\{\phi(V_{w}) : \phi$ is a chain of Whitehead
automorphisms of degree 0 throughout which the length of $V_{w}$
is constant, $|\phi(V_{w})|_{C_i}=|V_{w}|_{C_i} \ \text{for all} \
C_i$ with $C_i \neq C_1$, \ \text{and} \ $|\phi(V_{w})|_{C_1} \le
|\psi(V_{w})|_{C_1}$ for every $\psi$ which has the same property
as $\phi\}$; finally the same arguments as used in Claims ~6 and 7
in Case ~2.1 of Lemma ~3.2 are valid for $V_{w}$, since Hypothesis
~1.3 holds for $V_{w}$ if we only consider $C_i$'s of $\G_{v_{w}}$
such that $x_1 \notin C_i$ and $x_{n+1} \notin C_1$.

This observation allows us to apply the induction hypothesis on
$r$ to $\b \a_s \cdots \a_q (V_{w})$. Hence, there exist Whitehead
automorphisms $\g_1, \g_2, \dots, \g_h$ of $F_{n+3}$ such that
$$\b \a_s \cdots \a_q (V_{w}) = \g_h \cdots \g_2 \g_1 (V_{w}),
\tag 3.10
$$ where $r+1 \ge \deg \g_h \ge \deg \g_{h-1} \ge \cdots
\ge \deg \g_1$ (here note that there is no $\g_i$ of degree 1),
and $\sum\limits_{j=1}^m |\g_i \cdots \g_1
(v_j)|=\sum\limits_{j=1}^m |v_j|$ for all $i=1, \dots, h$.

As in Claim ~2, from each $\g_i$ we define a Whitehead
automorphism $\zeta_i$ of $F_n$. Let $k$ be such that $\deg
\zeta_j \le 1$ for $1 \le j < k$ and $\deg \zeta_j \ge 2$ for $k
\le j \le h$. Since $\b \a_s \cdots \a_q (V_{w})=V_{\s_{\ell+1}
\t_s \cdots \t_q (w)}$ and $\g_h \cdots \g_2 \g_1
(V_{w})=V_{\zeta_h \cdots \zeta_2 \zeta_1 (w)}$, we have by (3.10)
that
$$\s_{\ell+1} \t_s \cdots \t_q (w)= \zeta_h \cdots \zeta_2 \zeta_1 (w),$$
where $r+1 \ge \deg \zeta_h \ge \deg \zeta_{h-1} \ge \cdots \ge
\deg \zeta_k \ge 2$, and $|\zeta_i \cdots \zeta_1 (w)|=|w|$ for
$i=1, \dots, h$. Applying the base step for $r=1$ (that is, Lemma
~3.2) to $\zeta_{k-1} \cdots \zeta_1 \t_{q-1} \cdots \t_1 (u)$
completes the proof of Case ~3 . \qed

\heading 4. Proof of Theorem 1.5
\endheading
The aim of this section is to prove Theorem ~1.5. For a cyclic
word $w$ in $F_n$, let $N_k(w)$ denote the cardinality of the set
$\O_k(w)=\{\phi (w): \phi$ can be represented as a composition
$\t_s \cdots \t_1$ ($s \in \Bbb N$) of Whitehead automorphisms
$\t_i$ of $F_n$ of degree $k$ such that $|\t_i \cdots
\t_1(w)|=|w|$ for all $i=1, \dots, s \}$. Then bounding $N(u)$
reduces to bounding each $N_k(u)$, which is shown in the proof of
Theorem ~1.5 using the result of Theorem ~1.4. In Lemma ~4.1 we
bound $N_0(u)$. In Lemma ~4.2 we show that $N_k(u)$ for $k \ge 1$
is at most $N_0(V_u)$, where $V_u$ is a certain sequence of cyclic
words constructed from $u$, thus bounding $N_k(u)$ for $k \ge 1$.

\proclaim {Lemma 4.1} Let $u$ be a cyclic word in $F_n$. Then
$N_0(u)$ is bounded by a polynomial function of degree $n-2$ with
respect to $|u|$.
\endproclaim

\demo {Proof} Let $m_i$ be the number of occurrences of $x_i^{\pm
1}$ in $u$ for $i=1, \dots, n$. Clearly
$$N_0(u) \le N_0(x_1^{m_1}x_2^{m_ 2} \cdots x_n^{m_n}).$$
So it suffices to show that $N_0(x_1^{m_1}x_2^{m_ 2} \cdots
x_n^{m_n})$ is bounded by a polynomial function of degree $n-2$
with respect to $|u|$. For a cyclic word $v$ in $F_n$, define
$|v|_s$ as
$$|v|_s=\sum_{i=1}^n |v|_{C_i}.$$
Noting that $|x_1^{m_1}x_2^{m_ 2} \cdots x_n^{m_n}|_s=n$, put
$\Cal M=\{v: |v|_s=n$ and $v=\O_0(x_1^{m_1}x_2^{m_ 2} \cdots
x_n^{m_n})\}$, and $\Cal L= \{v: |v|_s>n$ and $v=\O_0
(x_1^{m_1}x_2^{m_ 2} \cdots x_n^{m_n})\}$. Obviously the
cardinality of $\Cal M$ is $(n-1)!$.

For the cardinality of $\Cal L$, let $v \in \Cal L$. Taking an
appropriate $u'\in \Cal M$ (note that $u'$ can be chosen as
follows: Write $v=x_{k_1} w_1 x_{k_2} w_2 \cdots x_{k_n} w_n$
(without cancellation), where $w_i$ is a (non-cyclic) word in
$\{x_{k_1}, \dots, x_{k_i}\}$; then
$u'=x_{k_1}^{m_{k_1}}x_{k_2}^{m_{k_2}} \cdots x_{k_n}^{m_{k_n}}$),
we have Whitehead automorphisms $\t_j=( A_j, a_j)$ of $F_n$ of
degree 0 such that
$$v=\t_s \cdots \t_1(u'),
\tag 4.1
$$
where $|\t_j \cdots \t_1(u')|=|u'|$ and $|\t_j \cdots \t_1(u')|_s
\ge |\t_{j-1} \cdots \t_1(u')|_s$ for all $j=1, \dots, s$. Then
for any $\t_i=(A_i, a_i)$ and $\t_j=(A_j, a_j)$ with $a_j \neq
a_i^{\pm 1}$, if we replace $\t_i$ and $\t_j$ by $(\bar A_i,
a_i^{-1})$ and $(\bar A_j, a_j^{-1})$, respectively, if necessary
so that $a_i^{\pm 1} \notin A_j$ and $a_j^{\pm 1} \notin A_i$,
then $A_i \cap A_j = \emptyset$. Hence by Case ~1.1.2 of Lemma
~3.1 that $\t_j \t_i \equiv \t_i \t_j$; thus (4.1) can be
re-written as
$$v=\t_{p t_p}^{q_{p t_p}} \cdots \t_{p1}^{q_{p 1}} \cdots \t_{1t_1}^{q_{1 t_1}} \cdots \t_{11}^{q_{11}}(u'),
\tag 4.2
$$
where $a_{k i}=a_{k i'}$ and $ A_{k i} \neq  A_{k i'}$ provided $i
\neq i'$; $a_{k' i} \neq a_{k i}^{\pm 1}$ and $(\t_{k'
t_{k'}}^{q_{k' t_{k'}}} \cdots \t_{k'1}^{q_{k'1}})(\t_{k
t_k}^{q_{k t_k}} \cdots \t_{k 1}^{q_{k 1}}) \equiv (\t_{k
t_k}^{q_{k t_k}} \cdots \t_{k 1}^{q_{k 1}}) (\t_{k' t_{k'}}^{q_{k'
t_{k'}}} \cdots \t_{k'1}^{q_{k'1}})$ provided $k \neq k'$. Here we
may assume by Case ~1.2.1 of Lemma ~3.1 that $ A_{k i} \subset
A_{k i'}$ if $i < i'$. Then $\t_{k i'}\t_{k i} \equiv \t_{k
i}\t_{k i'}$ by Case ~1.2.1 of Lemma ~3.1; hence $\t_{k' i'}
\t_{ki} \equiv \t_{ki}\t_{k'i'}$ for any $\t_{ki}$ and $\t_{k'i'}$
in chain (4.2).

\proclaim {Claim} The length of the chain of Whitehead
automorphisms on the right-hand side of (4.2) is at most $n-2$
without counting multiplicity, that is, $\sum \limits_{i=1}^p t_i
\le n-2$.
\endproclaim

\demo {Proof of the Claim} The proof proceeds by induction on the
number of subwords of $u'$ of the form $x_i^{m_i}$ which are fixed
throughout chain (4.2). For the base step, suppose that $u'$ has
two such subwords $x_{r_1}^{m_{r_1}}$ and $x_{r_2}^{m_{r_2}}$
(note that $u'$ must have at least two such subwords). The cyclic
word $u'$ can be written as $u'=x_{r_1}^{m_{r_1}} w$ (without
cancellation), where $w$ is a non-cyclic word that contains
$x_i^{m_i}$ for all $i \neq r_1$. Upon replacing $\t_{ij}$ by
$(\bar  A_{ij}, a_{ij}^{-1})$ if necessary, we may assume that
$x_{r_1}^{\pm 1} \notin  A_{ij}$ for all $\t_{ij}$ in chain (4.2).
Then the length of $w$ is constant throughout the chain and only
the subword $x_{r_2}^{m_{r_2}}$ of $w$ is fixed in passing from
$w$ to $\t_{p t_p}^{q_{p t_p}} \cdots \t_{p1}^{q_{p 1}} \cdots
\t_{1t_1}^{q_{1 t_1}} \cdots \t_{11}^{q_{11}}(w)$. It follows that
the length of this chain is precisely $(n-1)-1=n-2$ without
counting multiplicity. So the base step is done.

Now for the inductive step, suppose that $u'$ has $k$ subwords of
the form $x_i^{m_i}$ which are fixed throughout chain (4.2), say
$x_{r_1}^{m_{r_1}}, \dots, x_{r_k}^{m_{r_k}}$. Write the cyclic
word $u'$ as $u'=x_{r_1}^{m_{r_1}} w$ (without cancellation),
where $w$ is a non-cyclic word that contains $x_i^{m_i}$ for all
$i \neq r_1$. As above, upon replacing $\t_{ij}$ by $(\bar A_{ij},
a_{ij}^{-1})$ if necessary, we may assume that $x_{r_1}^{\pm 1}
\notin  A_{ij}$ for all $\t_{ij}$ in chain (4.2). We then have
that only the subwords $x_{r_2}^{m_{r_2}}, \dots,
x_{r_k}^{m_{r_k}}$ of $w$ are fixed in passing from $w$ to $\t_{p
t_p}^{q_{p t_p}} \cdots \t_{p1}^{q_{p 1}} \cdots \t_{1t_1}^{q_{1
t_1}} \cdots \t_{11}^{q_{11}}(w)$, where the length of $w$ is
constant throughout the chain.

Let $(w)$ be the cyclic word associated with $w$. If none of
$\t_{ij}$ in chain (4.2) is of the form either $(\Sig-x_{r_1}^{\pm
1}-x_g^{\pm 1}, x_g)$ or $(\Sig-x_{r_1}^{\pm 1}-x_g^{\pm 1},
x_g^{-1})$, then chain (4.2) can be applied to $(w)$ with $\t_{ij}
\neq 1$ on $(w)$ for every $\t_{ij}$ in the chain. Then by the
induction hypothesis applied to $(w)$, the length of the chain is
at most $(n-1)-2=n-3$ without counting multiplicity, as desired.
If one of $\t_{ij}$ in chain (4.2) is of the form either
$(\Sig-x_{r_1}^{\pm 1}-x_g^{\pm 1}, x_g)$ or $(\Sig-x_{r_1}^{\pm
1}-x_g^{\pm 1}, x_g^{-1})$, then we see that there can be only one
of $\t_{ij}$ of such a form, so that chain (4.2) can be applied to
$(w)$ with only one $\t_{ij}=1$ on $(w)$. This together with the
induction hypothesis applied to $(w)$ yields that the length of
chain (4.2) is at most $(n-1)-2+1=n-2$ without counting
multiplicity, as required. \qed
\enddemo

Obviously each multiplicity $q_{ij}$ is less than the number of
$a_{ij}^{\pm 1}$ occurring in $u$, so less than $|u|$. This
together with the Claim yields that the total number of chains of
Whitehead automorphisms with the same properties as in (4.2) is
less than $\binom r {n-2} |u|^{n-2}$, where $r$ is the number of
Whitehead automorphisms of $F_n$ of degree 0. Thus the cardinality
of $\Cal L$ is less than $(n-1)!\, \binom r {n-2}|u|^{n-2}$, and
therefore
$$N_0(x_1^{m_1}x_2^{m_ 2} \cdots x_n^{m_n}) =\# \Cal M +\# \Cal L \le (n-1)! + (n-1)!\, \binom r {n-2}|u|^{n-2},$$
which completes the proof the lemma.
\qed
\enddemo

\proclaim {Remark} The proof of Lemma ~4.1 can be applied without
further change if we replace consideration of a single cyclic word
$u$, the length $|u|$ of $u$, and the total number of occurrences
of $x_j^{\pm 1}$ in $u$ with consideration of a finite sequence
$(u_1, \dots, u_m)$ of cyclic words, the sum $\sum\limits_{i=1}^m
|u_i|$ of the lengths of $u_1, \dots, u_m$, and the total number
of occurrences of $x_j^{\pm 1}$ in $(u_1, \dots, u_m)$,
respectively.
\endproclaim

\proclaim {Lemma 4.2} Let $u$ be a cyclic word in $F_n$ that
satisfies Hypothesis 1.1. Then for each $k=1, \dots, n-1$,
$N_k(u)$ is bounded by a polynomial function of degree $n+3k-2$
with respect to $|u|$ (note that $k$ is at most $n-1$ by the
Remark after Lemma 2.2).
\endproclaim

\demo {Proof} Let $m_i$ be the number of occurrences of $x_i^{\pm
1}$ in $u$ for $i=1, \dots, n$, and let
$\ell_k=\sum\limits_{j=1}^k m_j$ for $k=1, \dots, n-1$. Write
$$u=y_1u_1y_2u_2 \cdots y_{\ell_k}u_{\ell_k} \hskip 0.2in \text {without cancellation},
\tag 4.3
$$
where for each $i=1, \dots, \ell_k$, $y_i= x_j$ or $y_i=x_j^{-1}$
for some $1 \le j \le k$, and $u_i$ is a (non-cyclic) subword in
$\{x_{k+1}, \dots, x_n\}^{\pm 1}$. Let $F_{n+3k}$ be the free
group on the set
$$\{x_1, \dots, x_n, x_{n+1}, \dots x_{n+k},
x_{2n+1}, \dots, x_{2n+k}, x_{3n+1}, \dots, x_{3n+k}\}.$$ From
(4.3) we construct a sequence $V_u=(v_1, \dots, v_{\ell_k})$ of
cyclic words $v_1, \dots, v_{\ell_k}$ in $F_{n+3k}$ with $\sum
\limits_{i=1}^{\ell_k} |v_i|=2|u|$ as follows: for each $i=1,
\dots, \ell_k$,
$$\aligned
&\text {if $y_i=x_j$ and $y_{i+1}=x_{j'}$, then $v_i=x_j u_i x_{3n+j'} u_i^{-1}$;}\\
&\text {if $y_i=x_j^{-1}$ and $y_{i+1}=x_{j'}$, then $v_i=x_{n+j} u_i x_{3n+j'} u_i^{-1}$;} \\
&\text {if $y_i=x_j$ and $y_{i+1}=x_{j'}^{-1}$, then $v_i=x_j u_i x_{2n+j'} u_i^{-1}$;} \\
&\text {if $y_i=x_j^{-1}$ and $y_{i+1}=x_{j'}^{-1}$, then
$v_i=x_{n+j} u_i x_{2n+j'} u_i^{-1}$,}
\endaligned
$$
where $y_{\ell_k+1}=y_1$.

\proclaim {Claim} For each Whitehead automorphism $\s$ of $F_n$ of
degree $k$ such that $|\s(u)|=|u|$, there exists a Whitehead
automorphism $\t$ of $F_{n+3k}$ of degree 0 such that
$\sum\limits_{i=1}^{\ell_k} |\t(v_i)|= \sum\limits_{i=1}^{\ell_k}
|v_i|$ and $\t(V_u)=V_{\s(u)}$.
\endproclaim

\demo {Proof of the Claim} Let $\s=(S, a)$ be a Whitehead
automorphism of $F_n$ of degree $k$ such that $|\s(u)|=|u|$. Upon
replacing $\s$ by $(\bar  S, a^{-1})$, we may assume that $\s=(
S,x_r)$. Note by Lemma ~2.2 that the index $r$ is bigger than $k$,
since $\deg \s=k$. Put $ S= T+ P+ Q$, where $ T= S \cap \{x_{k+1},
\dots, x_n\}^{\pm 1}$, $ P= S \cap \{x_1, \dots, x_k\}$ and $ Q= S
\cap \{x_1, \dots, x_k\}^{-1}$ (here note that $T=T^{-1}$, since
$\deg \s=k$).

Then we consider the Whitehead automorphism $\t=(T + P_1 + Q_1,
x_r)$ of $F_{n+3k}$ of degree ~0, where $P_1=\{x_i^{\pm 1},
x_{2n+i}^{\pm 1}| x_i \in P\}$ and $Q_1=\{x_{n+i}^{\pm 1},
x_{3n+i}^{\pm 1}| x_i^{-1} \in Q \}$. If the sequence $V_u=(v_1,
\dots, v_{\ell_k})$ of cyclic words $v_1, \dots, v_{\ell_k}$ in
$F_{n+3k}$ is constructed as above, then each newly introduced
letter $x_r^{\pm 1}$ in passing from $u$ to $\s(u)$ that remains
in $\s(u)$ produces two newly introduced letters $x_r^{\pm 1}$ in
passing from $V_u$ to $\t(V_u)$ that remain in $\t(V_u)$, and vice
versa. Also each letter $x_r^{\pm 1}$ in $u$ that is lost in
passing from $u$ to $\s(u)$ produces two letters $x_r^{\pm 1}$ in
$V_u$ that are lost in passing from $V_u$ to $\t(V_u)$, and vice
versa. This yields that $\sum\limits_{i=1}^{\ell_k} |\t(v_i)|=
\sum\limits_{i=1}^{\ell_k} |v_i|$. Moreover it is clear that
$\t(V_u)=V_{\s(u)}$. \qed
\enddemo

It is easy to see that if $u'\in \O_k(u)$ with $u' \neq u$, then
$V_{u'} \neq V_u$. This together with the Claim gives us that
$N_k(u) \le N_0((v_1, v_2, \dots, v_{\ell_k}))$. By the Remark
after Lemma ~4.1, $N_0((v_1, v_2, \dots, v_{\ell_k}))$ is bounded
by a polynomial function of degree $n+3k-2$ with respect to
$2|u|$, which completes the proof of the lemma. \qed
\enddemo

Finally we give a proof of Theorem ~1.5.

\proclaim {Proof of Theorem 1.5}
\endproclaim

Without loss of generality we may assume that $u$ was chosen from
the set $\{v \in \text {Orb}_{\text {Aut}F_n}(u): |v|=|u|\}$ so
that $u$ satisfies Hypothesis ~1.3. Let $v \in \text {Orb}_{\text
{Aut}F_n}(u)$ be such that $|v|=|u|$. By Whitehead's Theorem,
there exist Whitehead automorphisms $\pi$ of the first type and
$\s_1, \dots, \s_{\ell}$ of the second type such that $v=\pi
\s_{\ell} \cdots \s_1 (u)$, where $|\s_i \cdots \s_1 (u)|=|u|$ for
all $i=1, \dots, \ell$. Then by Theorem ~1.4, there exist
Whitehead automorphisms $\t_1, \dots, \t_s$ such that $v=\pi \t_s
\cdots \t_1 (u)$, where $n-1 \ge \deg \t_s \ge \deg \t_{s-1} \ge
\cdots \ge \deg \t_1$, and $|\t_j \cdots \t_1 (u)|=|u|$ for all
$j=1, \dots, s$ (here, note by the Remark after Lemma ~2.2 that
$\deg \t_s \le n-1$). This implies that
$$N(u) \le C N_0(u) N_1(u) \cdots  N_{n-1}(u),$$
where $C$ is the number of Whitehead automorphisms of the first
type of $F_n$ (which depends only on $n$). For each $k=0, 1,
\dots, n-1$, $N_k(u)$ is bounded by a polynomial function of
degree $n+3k-2$ with respect to $|u|$ by Lemmas 4.1 and 4.2.
Therefore, $N(u)$ is bounded by a polynomial function of degree
$n(5n-7)/2$ with respect to $|u|$, as required. \qed

\heading 5. Limitations
\endheading

We close this paper with a brief explanation why the presented
technique is incapable of covering the entire problem domain (e.g.
for $u=x_1^2x_2^2x_3^3x_4^4$ the presented arguments cannot be
applied). This amounts to explaining why condition (ii) of
Hypothesis ~1.1 cannot be dropped. As a matter of fact, in the
presented arguments, condition ~(ii) of Hypothesis ~1.1 played a
most essential role, without which all of our arguments except
Lemmas ~2.1 and 4.1 would have broke down. Owing to Lemma ~2.2
where we first used Hypothesis ~1.1 ~(ii), we were able to assume
throughout the paper that
$$j>i \ \text{when considering Whitehead automorphisms $(A, x_j^{\pm 1})$ of
degree $i$.} \tag 5.1
$$ This allowed us to exclude the worst case
such as $a \in B$, $a^{-1} \notin B$, $b \in A$ and $b^{-1} \notin
A$ in Lemma ~3.1, for which case there does not exist a
composition of Whitehead automorphisms of ascending degrees that
equals $(B, b)(A, a)$. Also we proceeded with the proofs of Lemma
~3.2 and Theorem ~1.4 based on (5.1). For instance, Claim ~1 in
the proof of Lemma ~3.2 yielded the existence $r$ such that
$a_r^{\pm 1} \in A_s \cap B$ in Case ~1.1, where we did not have
to worry about the case where $a_r \in A_s \cap B$ but $a_r^{-1}
\notin A_s \cap B$. Furthermore, the equality in the Claim in the
proof of Lemma ~4.2 would not have hold without (5.1).

\heading Acknowledgements
\endheading

The author would like to express her deep appreciation to
Professors S. V. Ivanov, I. Kapovich, A. G. Myasnikov, V.
Shpilrain for their comments and interests in this paper. The
author is also thankful to the referee for many valuable
suggestions which led to an improvement of this paper. The author
was partially supported by Research Institute for Basic Sciences,
Pusan National University (2005).

\heading References
\endheading

\roster

\item"1." P. J. Higgins and R. C. Lyndon, Equivalence of elements
under automorphisms of a free group, {\it J. London Math. Soc.}
{\bf 8} (1974), 254--258.

\item"2." I. Kapovich, P. E. Schupp and V. Shpilrain, Generic
properties of Whitehead's Algorithm and isomorphism rigidity of
random one-relator groups, {\it Pacific J. Math.} {\bf 223}
(2006), 113--140.

\item"3." B. Khan, The structure of automorphic conjugacy in the
free group of rank two, Computational and experimental group
theory, 115--196, {\it Contemp. Math.}, 349, Amer. Math. Soc.,
Providence, RI, 2004.

\item"4." R. C. Lyndon and P. E. Schupp, ``Combinatorial Group
Theory'', Springer-Verlag, New York/Berlin, 1977.

\item"5." J. McCool, A presentation for the automorphism group of
a free group of finite rank, {\it J. London Math. Soc.} {\bf 8}
(1974), 259--266.

\item"6." A. G. Myasnikov and V. Shpilrain, Automorphic orbits in free groups, {\it J. Algebra} {\bf 269} (2003), 18--27.

\item"7." J. H. C. Whitehead, Equivalent sets of elements in a
free group, {\it Ann. of Math.} {\bf 37} (1936), 782--800.

\endroster
\enddocument